[1994/12/01]
\documentclass{article}

\usepackage{arxiv}

\usepackage[utf8]{inputenc} 
\usepackage[T1]{fontenc}    
\usepackage{hyperref}       
\usepackage{url}            
\usepackage{booktabs}       
\usepackage{amsfonts}       
\usepackage{nicefrac}       
\usepackage{microtype}      
\usepackage{lipsum}

\usepackage{amsmath,amsfonts,amssymb,amsbsy,algorithm,algorithmic}
\usepackage{graphicx,footnote}
\usepackage{enumerate}
\usepackage{color}
\usepackage[displaymath,mathlines]{lineno}
\usepackage{natbib}
\usepackage{etoolbox}

\title{Enhanced parallelization of the incremental 4D-Var data assimilation algorithm using the Randomized Incremental Optimal Technique (RIOT)}

\author{
  Nicolas Bousserez \thanks{Now at ECMWF, Reading, UK} \\
  University of Colorado \\
  Department of Mechanical Engineering \\
  Boulder, CO, USA \\
  \texttt{nicolas.bousserez@ecmwf.int} 
   \And
Jonathan J. Guerrette \thanks{ Now at NCAR, Mesoscale and Microscale Meteorology Laboratory, Boulder, CO, USA} \\
 University of Colorado \\
  Department of Mechanical Engineering \\
  Boulder, CO, USA \\
NOAA, Chemical Sciences Division \\
 Boulder, CO, USA \\
    \And
 Daven K. Henze \\
 University of Colorado \\
 Department of Mechanical Engineering\\
Boulder, CO, USA }
\chardef\bslash=`\\ 

\newcommand{\st}{\sigma}

\newcommand{\eval}[2][\right]{\relax
  \ifx#1\right\relax \left.\fi#2#1\rvert}

\newcommand{\algorithmicdoinparallel}{\textbf{do in parallel}}
\makeatletter
\AtBeginEnvironment{algorithmic}{%
  \newcommand{\FORALLP}[2][default]{\ALC@it\algorithmicforall\ #2\ %
    \algorithmicdoinparallel\ALC@com{#1}\begin{ALC@for}}%
}

\begin{document}

\markboth{RIOT 4D-Var}
 {Enhanced parallelization of the incremental 4D-Var data assimilation algorithm using the Randomized Incremental Optimal Technique (RIOT)}

\renewcommand{\sectionmark}[1]{}

\begin{abstract}
\normalsize{Incremental 4D-Var is a data assimilation algorithm used routinely at operational numerical weather predictions centers worldwide. The algorithm solves a series of quadratic minimization problems (inner-loops) obtained from linear approximations of the forward model around non-linear trajectories (outer-loops). Since most of the computational burden is associated with the inner-loops, many studies have focused on developing computationally efficient algorithms to solve the least-square quadratic minimization problem, in particular through time parallelization. This paper implements a new method for parallelizing incremental 4D-Var, the Randomized Incremental Optimal Technique (RIOT), which replaces the traditional sequential conjugate gradient (CG) iterations in the inner-loop of the minimization with fully parallel randomized singular value decomposition (RSVD) of the preconditioned Hessian of the cost function. RIOT is tested using the standard Lorenz-96 model (L-96) as well as two realistic high-dimensional atmospheric source inversion problems based on aircraft observations of black carbon concentrations. A new outer-loop preconditioning technique tailored to RSVD was introduced to improve convergence stability and performance. Results obtained with the L-96 system show that the performance improvement from RIOT compared to standard CG algorithms increases significantly with non-linearities. Overall, in the realistic black carbon source inversion experiments, RIOT reduces the wall-time of the 4D-Var minimization by a factor 2-3, at the cost of a factor 4-10 increase in energy cost due to the large number of parallel cores used. Furthermore, RIOT enables reduction of the wall-time computation of the analysis error covariance matrix by a factor 40 compared to a standard iterative Lanczos approach. Finally, as evidenced in this study, implementation of RIOT in an operational numerical weather prediction system will require a better understanding of its convergence properties as a function of the Hessian characteristics and, in particular, the degree of freedom for signal (DOFs) of the inverse problem. }
\end{abstract}

\maketitle

\section{Introduction}
Inverse problems and data assimilation (DA) in  geophysical sciences are often challenging owing to the large size of the parameter and observation spaces. In the context of non-linear Bayesian inference, an analytical solution is often not computationally tractable.  Techniques thus have to be designed to efficiently approximate the solution to the inverse problem, which leads to a proliferation of approaches as described in overviews of methods applied throughout geophysics \citep{Carrassi2018, Asch2016, tarantola2005inverse}, or those specifically related to numerical weather prediction \citep{kalnay2003atmospheric}, remote sensing \citep{Rodgers2000}, atmospheric chemistry \citep{Bocquet2015} and biogeochemistry\citep{Rayner2018}.  The subset of techniques known as variational (Var) approaches seek to solve for the maximum a posteriori solution (MAP) of the Bayesian problem using iterative minimization procedures. 
%
%
%
%
%
In particular, one such widely used DA technique is the incremental 4D variational approach \citep{courtier1994strategy}. This method entails a series of linearizations of the system to formulate quadratic optimization problems, whose solutions (the increments) are often sought using conjugate gradient (CG) methods.  Each quadratic minimization, called an ``inner-loop'', is followed by an ``outer loop'' iteration wherein the full non-linear model is integrated to propagate the updated increment.  The incremental 4D-Var approach is computationally efficient in that the inner-loop minimizations may be performed using simplified, dimension-reduced, and/or coarser resolution tangent-linear (TL) and adjoint (AD) models \citep{Lawless2008}.
%
%
%
It has thus been used at operational numerical weather prediction (NWP) centers such as the UK Met Office \citep{Rawlins2007},  the Meteorological Service of Canada \citep{Laroche2005}, M\'{e}t\'{e}o-France \citep{Gauthier2001}, the Japan Meteorological Agency \citep{Honda2005}, and the European Centre for Medium-Range Weather Forecasts \citep{Rabier2000}.

 Despite the widespread use of the incremental 4D-Var approach, a drawback (especially when compared to ensemble-based DA techniques \citep[e.g.,][]{Evensen2009,Houtemaker2016}, 
 %
 is that this method does not easily lend itself to parallelization (i.e., beyond using parallel computing to speed up the evaluation of any single  model evaluation).   Spatial parallelization of variational DA techniques has long been considered \citep[e.g.,][]{Elbern99,Rantakokko97,Tremolet96}.  Moreover, methods have  been developed to exploit time parallelism in 4D-Var through sub-division of the assimilation time window \citep{Rao2016} or through the saddle-point formulation in the weak-constraint algorithm \citep{fisher2017}, although the convergence of the latter cannot always be guaranteed \citep{Gratton2018}. Most recently, \citet{Mercier18,Mercier19} have used block-krylov methods to accelerate an ensemble of variational DA in the context of 3D-Var, where all observations are assimilated at the same time $t_0$ (i.e., the 4D model dynamic is not considered in the assimilation). In this approach, Krylov subspaces computed in parallel for different perturbed 3D-Var ensemble members are concatenated to form a larger Krylov subspace which is used to minimize simultaneously all members, resulting in a speedup of the optimization of 13-45\%. However, this framework has not been applied to an ensemble of 4D-Var minimizations, where larger non-linearities are expected that could limit the efficiency of the approach due to significant variation of the Hessian across ensemble members.

In  this study, we evaluate a different and novel approach to parallelizing strong-constrained incremental 4D-Var algorithms, which was recently proposed by \citet{bous18} -- the Randomized Incremental Optimal Technique (RIOT). 
The key component of the algorithm is the use of randomized singular vector decomposition (RSVD) methods \citep{Halko11}, which enables parallel model integrations in place of the traditional sequential minimization steps of the inner loops. 
RIOT thus introduces a new level of parallelization for variational DA that is advantageous when the available  computing resources exceed that which can be efficiently used by any single forward model evaluation (e.g., through parallel computation with MPI communications).

Here we evaluate the performance of RIOT against standard incremental 4D-Var techniques that perform the inner-loop minimizations using iterative CG algorithms. Our numerical experiments are based on a Lorenz-96 (hereafter L-96) system, as well as a high-dimensional source inversion problem  (recovering emissions of black carbon aerosol based on ambient concentration measurements). In Section \ref{theory}, we present the theoretical background of the incremental 4D-Var method and introduce the different algorithms evaluated in this study.  In Section \ref{num_exp}, we present and discuss results from numerical DA experiments and assess the potential of RIOT by analyzing its computational and convergence characteristics compared  with those of traditional iterative minimization methods. Finally, a summary of our results and a discussion on the limitations and potential of RIOT for use in operational  4D-Var DA systems is provided in the conclusion.

\section{Theory and Methods}
\subsection{Incremental 4D-Var}
\label{theory}
In data assimilation,  we seek the analysis, $\mathbf{x}_a$, that is the Maximum a Posteriori (MAP) solution of a Bayesian problem, whose  posterior distribution is proportional to the product of the prior and likelihood probability distribution functions (pdfs). It is common to assume that the prior and likelihood pdfs are  Gaussian. More specifically, we assume that
 $\mathbf{x} \sim \mathcal{N}(\mathbf{x}^b,\mathbf{B})$ is the vector of parameters with prior distribution associated with the background $\mathbf{x}^b$ and  covariance error matrix $\mathbf{B}\in\mathbb{R}^{n \times n}$, and that $\mathbf{y}$ is a vector of observations distributed in time and space (4D) with likelihood function $\mathbf{y}|\mathbf{x} \sim \mathcal{N}(H(\mathbf{x}),\mathbf{R})$ associated with a forward model $H$ and an error covariance matrix $\mathbf{R}\in\mathbb{R}^{r \times r}$. 
 Here we use the term parameters in a general sense that includes both model state prognostic variables and non-state variables (e.g., boundary conditions). The analysis $\mathbf{x}_a$ can be calculated by minimizing the following cost function:
\begin{eqnarray}
 \label{eq2}
J(\mathbf{x}) =&\frac{1}{2}(H(\mathbf{x})-\mathbf{y})^T\mathbf{R}^{-1}(H(\mathbf{x})-\mathbf{y}) \\  \nonumber
&+\frac{1}{2}(\mathbf{x}-\mathbf{x}_{b})^T\mathbf{B}^{-1}(\mathbf{x}-\mathbf{x}_{b}).
\end{eqnarray}
When the parameter and observation vectors are both high-dimensional and the integration of $H$ is computationally expensive, the minimum of (\ref{eq2}) can be obtained using an iterative algorithm, which is the aforementioned 4D-Var approach. Two popular gradient-based optimization algorithms for large-scale non-quadratic problems are the limited-memory Broyden-Fletcher-Goldfarb-Shanno (L-BFGS) method \citep{Liu89,byrd1995limited}
%
%
%
and the truncated Gauss-Newton method.  

Incremental 4D-Var \citep{courtier1994strategy} is a special case of the latter \citep{Lawless2005,Gratton2007}
%
that approximately \citep{Lawless2008} minimizes a sequence of quadratic cost functions to solve the non-quadratic problem. In short, each quadratic minimization (inner-loop) solution provides a so-called increment that is added to the current parameter vector, from which another non-linear integration is performed (outer-loop) to update the model trajectory. A linearization of the model around that new model trajectory leads to the next quadratic cost function minimization. The process continues until a convergence criteria is reached. The quadratic minimization in incremental 4D-Var consists of solving the following problem: 
 \begin{align}
 \label{eq3}
&\delta \mathbf{z}_a = \min_\mathbf{\delta z} \tilde{J}(\mathbf{\delta z}) =\frac{1}{2}(\mathbf{d}-\mathbf{H}\mathbf{L\delta z})^T\mathbf{R}^{-1}(\mathbf{d}-\mathbf{H}\mathbf{L\delta z})+\frac{1}{2}\delta\mathbf{z}^T\delta\mathbf{z} \\
& \delta\mathbf{x}_a = \mathbf{L} \delta \mathbf{z}_a \\
& \mathbf{B}=\mathbf{LL}^T  \\
  &\mathbf{d}=\mathbf{y}-H(\mathbf{x}),
\end{align}
where $\delta\mathbf{x}_a$ and $\mathbf{d}$ are the increment and innovation vectors, respectively, and $\mathbf{H}$ denotes the Jacobian evaluated at $\mathbf{x}$. The minimization problem (\ref{eq3}) is usually solved using an iterative CG algorithm \citep{Chao92}, for which only the product of the Jacobian or its transpose with a vector is required. Therefore, only tangent linear (TL) and adjoint (AD) model integrations are necessary to solve the minimization problem (\ref{eq3}).

\subsection{Principle of the Methods}
\label{methods}

\paragraph{Approximating the Inverse Hessian}
It has been shown that the CG method is strongly connected to the Lanczos algorithm \citep{Golub89}. In particular, the solution of (\ref{eq3}) obtained from $q$ iterations of a CG minimization can be written as follow:
 \begin{align}
 \label{eq4}
 \delta \mathbf{z}^a_q=\left [\mathbf{Q} \mathbf{Q}^T ( \mathbf{A}+\mathbf{I} )\mathbf{Q} \mathbf{Q}^T \right ]^{+}\mathbf{b},
  \end{align}
where $\mathbf{A}$ is the prior-preconditioned Hessian $=\mathbf{L}^T \mathbf{H}^T \mathbf{R}^{-1}\mathbf{H}\mathbf{L}$, $\mathbf{Q}$ is an orthonormal basis of $\mathbf{Y}=\{\mathbf{A}^q\mathbf{b},...,\mathbf{A} \mathbf{b} ,\mathbf{b}\}$, $^+$ denotes the pseudo-inverse and $\mathbf{b}=\mathbf{L}^T\mathbf{H}^T \mathbf{R}^{-1}\mathbf{d}$.

Furthermore, the pseudo-inverse of a matrix can be expressed as a function of its singular decomposition, i.e.:
\begin{align}
  \label{eq5}
  \mathbf{M}^+=\sum_{i=1}^p \mu_i^{-1} \mathbf{v}_i \mathbf{w}_i^T,  \, \mathbf{M} \in \mathbb{R}^{n \times m},\, p\le \min(m,n)
  \end{align}
where $(\mathbf{v}_i,\mathbf{w}_i)$ and $\mu_i$ represent the singular vectors and singular values of $\mathbf{M}$, respectively. Therefore, one sees from (\ref{eq4}) and (\ref{eq5}) that $q$ iterations of a CG minimization involves an implicit approximated eigenvalue decomposition of $\mathbf{A}$. 

In \citet{bous18}, following earlier work by \citet{BuiThanh2012}, an alternative to using the CG method for minimizing (\ref{eq3}) was proposed, in which the Krylov subspace $\mathbf{Y}$ approximating the range of $\mathbf{A}$ is replaced by a randomized range finder $\mathbf{Y}=\{\mathbf{A}\boldsymbol{\epsilon}_q,...,\mathbf{A} \boldsymbol{\epsilon}_1 \}$, where the $\boldsymbol{\epsilon}_i \in \mathbb{R}^{n} $ are random vectors drawn from a normal distribution \citep{Halko11}. The advantage of this approach is that it presents a new dimension for parallelization of the algorithm, since each element of $\mathbf{Y}$ can be computed independently. Another feature of RIOT is that it uses an adaptive increment at each inner-loop based on the approximated spectra of $\bf A$ that ensures statistical optimality. The interested reader can refer to \citet{bous18} for more details on the optimal adaptive approach. 

In addition to these purely sequential and parallel methods, others such as the randomized power iteration (e.g., \citet{Halko11}) and block-Krylov approaches use both iteration and matrix-block multiplication.  \citet{Musco2015} show that their version of a block-Krylov approach is more efficient for matrix decomposition than the randomized power iteration, which is why we explore the block-Lanczos algorithm as an alternative to CG or RSVD.  Block-Lanczos \citep{Golub89} is designed to solve a linear-least-squares problem where there is statistical uncertainty in the right-hand-side of the equation.  In the sequence of quadratic minimization problems of incremental 4D-Var, the right-hand-side is the cost function gradient, the uncertainty of which can be approximated by realizations of the observation ($\mathbf{y}$) and prior ($\mathbf{x}_b$). In that case the range approximation of $\mathbf{A}$ becomes $\{\mathbf{Y}=\mathbf{Y}_k,...,\mathbf{Y}_1\}$, where $\mathbf{Y}_i=\{\mathbf{A}^q\mathbf{b}_i,...,\mathbf{A} \mathbf{b}_i ,\mathbf{b}_i\}$ and $\mathbf{b}_i$ is an independent realization of the gradient (see Alg. \ref{alg:Block_lanczos}). This approach thus combines aspects of both  RSVD (parallelization of the prior-preconditioned Hessian-vector products) and  CG methods (construction of Krylov subspaces). Its advantages are that it usually requires a smaller number of iterations than CG and fewer ensemble members than RSVD, therefore enabling the algorithm to adjust to available computational resources.

\paragraph{Preconditioning}
The quadratic minimization (\ref{eq3}) uses the change of variable $\mathbf{x}=\mathbf{L}\delta\mathbf{z}$. This approach allows one to: 1) implicitly invert the background covariance matrix $\mathbf{B}$ in (\ref{eq2}), since this large matrix cannot be represented explicitly; 2) improve the convergence rate of the minimization by reducing the condition number of the Hessian of the cost function. Reducing the condition number of the Hessian through a change of variable is called preconditioning the problem. An additional level of preconditioning is possible by exploiting the approximations of the Hessian obtained as a by-product in each inner-loop minimization \citep{Tshimanga2008}. More specifically, the preconditioned problem consists of solving: 
\begin{align}
  \label{eq6}
 \mathbf{P}^T( \mathbf{A}+\mathbf{I} )\mathbf{P} \delta \mathbf{z}_a= \mathbf{P}^T\mathbf{b} \\
 \kappa( \mathbf{P}^T \mathbf{A}  \mathbf{P}) < \kappa(\mathbf{P}) ,
  \end{align}
 where $\kappa(.)$ represents the condition number of a matrix, i.e., the ratio between its largest and smallest eigenvalues. Ideally, $\mathbf{P}^T \mathbf{A}  \mathbf{P}=\mathbf{I}$, in which case the CG minimization converges in one iteration. This is because by construction the CG algorithm has converged when the Krylov subspace $\mathbf{Y}$ becomes invariant under multiplication by $\mathbf{A}$. Therefore, a good strategy is to choose a preconditioner $\mathbf{P}$ such that: $\mathbf{PP}^T\approx \mathbf{A}^+$. A standard approach is to use the following spectral preconditioner:
 \begin{align}
  \label{eq7}
 \mathbf{P}=\mathbf{I}+\sum_{i=1}^q (\lambda_i^{-1/2}-1) \mathbf{v}_i\mathbf{v}_i^T,
  \end{align}
  where the $\mathbf{v}_i$ and $\lambda_i$ approximate the eigenvectors and eigenvalues of $\mathbf{A}$, respectively \citep{Tshimanga2008}. The efficiency of the spectral preconditioner has been verified, e.g., in the ECMWF 4D-Var operational DA system, where it decreases the wall-time of the minimization by about 20\% compared to the non-preconditioned case (personal communication).
  
  On the other hand, the RIOT algorithm relies on an approximation of the inverse Hessian (see Eq. \ref{eq4}) using random samples of the range of $\mathbf{A}$. As a consequence, the subspace (defined by $\mathbf{Q}$) in which the cost function is minimized is agnostic to the right-hand term of (\ref{eq6}) and will thus be suboptimal (i.e., for a given rank the residual error is not minimized). Furthermore, assuming weak non-linearities, the inverse Hessian approximations computed from RIOT in (\ref{eq4}) will tend to be similar across outer iterations, although the modes resolved in previous outer iterations will have only small contributions to the current inner loop solution. This may prevent the algorithm from further reducing the cost function beyond the first outer iteration and thus, convergence of the minimization. 
  
  An additional concern is related to the eigenvalue spectra.  The spectral preconditioner (\ref{eq6}) transformation tends to cluster around 1 the eigenvalues associated with the modes resolved in previous outer iterations, which can be used to filter them out in the  RSVD Hessian approximation of subsequent iterations to more efficiently (using fewer ensemble members) extract the remaining unresolved modes with eigenvalues greater than 1. 
  However,  the accuracy of the RSVD method itself is somewhat at odds with this approach.  While RSVD best approximates the dominant eigenmodes of the targeted matrix \citep{Halko11}, the accuracy of the estimation is degraded for flat eigenvalue spectra, i.e., when the rate of decrease of the eigenvalues is small. As a consequence, the decrease in the Hessian condition number across outer iterations from the spectral preconditioner can result in poor minimization performances in RIOT. One solution to mitigate the interference of previously resolved modes in the RSVD approximation is to restrict the sampling of the Hessian matrix $\mathbf{Y}=\{\mathbf{A}\boldsymbol{\epsilon}_q,...,\mathbf{A} \boldsymbol{\epsilon}_1 \}$ to the subspace orthogonal to those modes. This improvement will be referred to as the rotation method in the following. The rotation technique is presented in Alg. \ref{alg:riot_primal_precond}.

\paragraph{The Adaptive Approach}
The RIOT algorithm computes a low-rank RSVD approximation of the (prior-preconditioned) Hessian, from which an approximated posterior increment and posterior error covariance matrix are derived.  Based on the approximated eigendecomposition of the Hessian, two updates are possible for the posterior error covariance matrix, as follows \citep{bous18}:
\begin{align}
\label{eq8}
\mathbf{P}_a& \approx \mathbf{B}-\mathbf{L}^T\left (  \sum_{i=1}^m \lambda_i(1+\lambda_i)^{-1}\mathbf{v}_i\mathbf{v}_i^T) \right ) \mathbf{L} &\text{  (LRU)} \\
\label{eq9}
\mathbf{P}_a& \approx \mathbf{L}^T\left (  \sum_{i=1}^m (1+\lambda_i)^{-1}\mathbf{v}_i\mathbf{v}_i^T) \right ) \mathbf{L}, &\text{  (LRA)}
\end{align}
where the $(\mathbf{v}_i,\lambda_i,\, i=1,...,m)$ are the eigenvectors and values, respectively, of the prior-preconditioned Hessian $\mathbf{A}=\mathbf{L}^T \mathbf{H}^T \mathbf{R}^{-1}\mathbf{H}\mathbf{L}$, in decreasing eigenvalue order, and where $m$ defines the truncation number for the eigenvalue decomposition. Here the updates LRU and LRA stand for low-rank update and low-rank approximation, respectively, since (\ref{eq8}) consists in an update to the prior covariance matrix $\mathbf{B}$ and is therefore full-rank, while (\ref{eq9}) is a rank-$m$ approximation to the posterior covariance matrix. The corresponding LRU and LRA posterior increment approximations are obtained by using (\ref{eq8}) and (\ref{eq9}), respectively, in the increment formula:
\begin{align}
\label{eq10}
\delta \mathbf{x}_a = \mathbf{L}\delta \mathbf{z}_a=\mathbf{P}_a\mathbf{b}
\end{align}

Assuming an exact truncated eigenvalue decomposition (TSVD) of the Hessian is used in  (\ref{eq8}) and (\ref{eq9}), statistical optimality results can be used to design an adaptive approach to select  either the LRU or LRA updates based on the approximated eigenvalue spectra. Thus, the following selection strategy is used in the RIOT algorithm: 
\begin{align}
\label{eq11} 
&\textbf{For } \lambda_k\ge1,\textbf{ use LRA posterior updates}  \\ 
&\textbf{For }  \lambda_k<1,\textbf{ use LRU posterior updates}  
\end{align}

The reader can refer to Bousserez and Henze (2018) for more information on the theory.

\subsection{Algorithms}
We recall below the RIOT algorithm, which we shall evaluate in this study. A variant of RIOT for the dual formulation of 4D-Var (RIOT-PSAS) is provided in the Appendix. The standard incremental CG algorithm (VarCG), as well as the Block-Lanczos method (VarBL) are also reformulated below so as to facilitate direct comparisons between methods. VarCG is presented using the equivalent Lanczos recurrence in the inner loop for simple comparison to the block variant and to recall the derivation of the approximate Hessian eigen modes. The need for re-orthogonalization is well-documented for both sequential and block Krylov methods \citep{Golub89}, and these are applied for both VarCG and VarBL. For the sake of readability, the spectral preconditioning steps are not included for the VarCG, VarBL and RIOT-PSAS methods, but only presented in the RIOT algorithms.
The VarCG minimization algorithm can be written as follows:
\begin{algorithm}
  \caption{One Cycle of Incremental 4D-Var with Conjugate-Gradient Minimization (VarCG)}\label{alg:lanczos}
    \begin{algorithmic}[1]                    
      \STATE Start with $\mathbf{x}_{0}=\mathbf{x}_b$, $\mathbf{v}_0=\mathbf{0}$, $\mathbf{B}=\mathbf{L}_\mathbf{B}\mathbf{L}_\mathbf{B}^\top$.
      \STATE Choose $m$ (number of iterations).
      \STATE Choose $k_f$ (number of outer loops).
      \FOR{$k=1, 2,\hdots, k_f$}         
            \STATE Integrate and store trajectory $H\left(\mathbf{x}_{k-1}\right)$.                 
            \STATE Compute and store $\mathbf{d}_{k-1}=H(\mathbf{x}_{k-1})-\mathbf{y}$.
            \STATE Let $\hat{\mathbf{A}} \equiv \mathbf{L}_\mathbf{B}^\top\mathbf{H}_{k-1}^\top\mathbf{R}^{-1}\mathbf{H}_{k-1}\mathbf{L}_\mathbf{B} \in \mathbb{R}^{n \times n}$, where $\mathbf{H}_{k-1}$ is the tangent-linear and $\mathbf{H}_{k-1}^\top$ the adjoint model at $\mathbf{x}_{k-1}$.
            \STATE Compute the gradient $\mathbf{b}\equiv\mathbf{L}_\mathbf{B}^\top\mathbf{H}_{k-1}^\top\mathbf{R}^{-1}\mathbf{d}_{k-1}+\mathbf{v}_{k-1}$.
            \STATE Set $\mathbf{r}_{0} = \mathbf{b}$.
            \STATE Set $\beta_{0} = \|\mathbf{r}_{1}\|_2$.
            \STATE $\mathbf{q}_0 = \mathbf{0}$
            \FORALL {$i\in\{1,\hdots,m\}$}
                  \STATE Set $\mathbf{q}_i = \mathbf{r}_{i-1} / \|\beta_{i-1}\|$.
                  \STATE Set $\mathbf{w} = \hat{\mathbf{A}}\mathbf{q}_{i}$.
                  \STATE Set $\alpha_i = \mathbf{q}_i^\top \mathbf{w}$.
                  \STATE Set $\mathbf{r}_{i} = \mathbf{w} - \alpha_i \mathbf{q_i} - \beta_i \mathbf{q}_{i-1}$.
                  \STATE Othogonalize $\mathbf{r}_{i}$ with respect to $\mathbf{r}_{0:i-1}$.
                  \STATE Set $\beta_{i} = \|\mathbf{r}_{i}\|_2$.
            \ENDFOR
            \STATE Form the tridiagonal matrix $\mathbf{T}\in\mathbb{R}^{m \times m}$ with $\alpha_{1 : m}$ on the diagonal and $\beta_{1 : m-1}$ on the subdiagonal and superdiagonal.
            \STATE Perform eigendecomposition $\mathbf{T} = \mathbf{Z}\mathbf{\Lambda}\mathbf{Z}^\top$.
            \STATE Compute  $\mathbf{U}=\mathbf{Q}\mathbf{Z}$, where $\mathbf{Q}=[\mathbf{q}_1 \hdots \mathbf{q}_m]$ is the Krylov subspace.
            \STATE $\delta\mathbf{v}= -  \left (\sum_{i=1}^m \frac{1}{1+\lambda_i}\mathbf{u}_i\mathbf{u}_i^\top \right ) \mathbf{b} $
            \STATE $\mathbf{v}_{k} = \mathbf{v}_{k-1}+\delta\mathbf{v}$
            \STATE $\mathbf{x}_{k} = \mathbf{x}_{k-1}+\mathbf{L}_\mathbf{B}\delta\mathbf{v}$
      \ENDFOR
    \end{algorithmic}
%
\end{algorithm}

\newpage
As explained previously, the VarBL method realizes a similar algorithmic parallelization as RIOT, but also requires iterative constructions of Krylov subspaces similar to VarCG. The inner-loop minimization follows the Block Lanczos algorithm given by \citet{Golub89}, where the starting block is comprised of realizations of the cost function gradient.  VarBL can be written as follows:

 \begin{algorithm}
  \caption{One Cycle of Incremental 4D-Var with Block-Lanczos Minimization (VarBL)}\label{alg:Block_lanczos}
    \begin{algorithmic}[1]                    
        \STATE Start with $\mathbf{x}_{0}=\mathbf{x}_b$, $\mathbf{v}_0=\mathbf{0}$, $\mathbf{B}=\mathbf{L}_\mathbf{B}\mathbf{L}_\mathbf{B}^\top$, $\mathbf{R}=\mathbf{L}_\mathbf{R}\mathbf{L}_\mathbf{R}^\top$.
        \STATE Choose $m$ (number of iterations).
        \STATE Choose $l$ (number of gradient samples).
        \STATE Choose $k_f$ (number of outer loops).
        \FOR{$k=1, 2,\hdots, k_f$}
                  
              \STATE Draw $\epsilon_{j}\in\mathbb{R}^{p \times l}\sim\mathcal{N}(0,1)$. 
              \STATE Draw $\epsilon_{b,j}\in\mathbb{R}^{n \times l}\sim\{\mathcal{N}(0,1), j=[2,\ldots,l]\}$.
              \STATE Set $\epsilon_{b,1}=\mathbf{0}$.
              \FORALLP {$j\in\{1,\hdots,l\}$}
                    \STATE Integrate and store trajectory $H\left(\mathbf{x}_{k-1} + \mathbf{L}_\mathbf{B}\epsilon_{b,j}\right)$.
                    \STATE Compute $\mathbf{d}_{{k-1},j}=H(\mathbf{x}_{k-1} + \mathbf{L}_\mathbf{B}\epsilon_{b,j})-\mathbf{y}+\mathbf{L}_\mathbf{R}\epsilon_{j}$.
                    \STATE Compute the gradient $\mathbf{b}_j\equiv\mathbf{L}_\mathbf{B}^\top\mathbf{H}_{k-1,j}^\top\mathbf{R}^{-1}\mathbf{d}_{{k-1},j} + \mathbf{v}_{k-1} + \epsilon_{b,j}$, where $\mathbf{H}_{k-1,j}^\top$ is the adjoint model evaluated at $\mathbf{x}_{k-1} + \mathbf{L}_\mathbf{B}\epsilon_{b,j}$
              \ENDFOR
              \STATE Let $\hat{\mathbf{A}} \equiv \mathbf{L}_\mathbf{B}^\top\mathbf{H}_{k-1}^\top\mathbf{R}^{-1}\mathbf{H}_{k-1}\mathbf{L}_\mathbf{B} \in \mathbb{R}^{n \times n}$, where $\mathbf{H}_{k-1}$ and $\mathbf{H}_{k-1}^\top$ are the tangent-linear and adjoint models evaluated at $\mathbf{x}_{k-1}$.
              \STATE Establish $\mathbf{Q}[ : , 1:l ]$ with the ortho-normalized gradient vectors $\mathbf{b}_j, j=[1,\ldots,l]\}$
              \FORALL {$i\in\{1,\hdots,m\}$}
                    \FORALLP {$j\in\{1,\hdots,l\}$}
                          \STATE $\mathbf{Y}[ : , j + (i-1)*l] = \hat{\mathbf{A}}\mathbf{Q}[: , j + (i-1)*l]$                  
                    \ENDFOR
                    \STATE $\mathbf{M}_{i} = \mathbf{Q}[: , 1 + (i-1)*l : i*l ]^\top \mathbf{Y}[ : , 1 + (i-1)*l : i*l]$
                    \IF { $i < m$ }                  
                          \STATE Calculate 2-term residual $\mathbf{r} = \mathbf{Y}[ : , 1 + (i-1)*l : i*l] - \mathbf{Q}[: , 1 + (i-1)*l : i*l] \mathbf{M}_{i}$
                          \IF { $i > 1$ }
                                \STATE Complete Block Lanczos 3-term recurrence: $\mathbf{r} = \mathbf{r} - \mathbf{Q}[: , 1 + (i-2)*l : (i-1)*l] \mathbf{W}_{i-1}^\top$
                          \ENDIF
                          \STATE Evaluate QR decomposition, $\mathbf{Q}[: , 1 + i*l : (i+1)*l] \mathbf{W}_{i} = \mathbf{r}$
                          \STATE Ortho-normalize $\mathbf{Q}[: , 1 + i*l : (i+1)*l]$ relative to $\mathbf{Q}[: , 1 : i*l]$ 
                    \ENDIF
              \ENDFOR
              \STATE Form block-tridiagonal matrix $\mathbf{T}\in\mathbb{R}^{ml \times ml}$ with $\mathbf{M}_{1 : m}$ on the diagonal and $\mathbf{W}_{1 : m-1}$ and ${\mathbf{W}^\top}_{1 : m-1}$ on the subdiagonal and superdiagonal, respectively.
              
              \STATE Perform eigendecomposition $\mathbf{T} = \mathbf{Z}\mathbf{\Lambda}\mathbf{Z}^\top$.
              \STATE Compute  $\mathbf{U}=\mathbf{Q}\mathbf{Z}$.
              \STATE $\delta\mathbf{v}= -  \left (\sum_{i=1}^m \frac{1}{1+\lambda_i}\mathbf{u}_i\mathbf{u}_i^\top \right ) \mathbf{b} $
              \STATE $\mathbf{v}_{k} = \mathbf{v}_{k-1}+\delta\mathbf{v}$
              \STATE $\mathbf{x}_{k} = \mathbf{x}_{k-1}+\mathbf{L}_\mathbf{B}\delta\mathbf{v}$
        \ENDFOR
    \end{algorithmic}
%
\end{algorithm}

\newpage
 Finally, the RIOT algorithm, which includes preconditioning and sampling improvement (rotation) options, is as follows:
  \begin{algorithm}
  \caption{  One Cycle of the Randomized Incremental Optimal Technique (RIOT)}\label{alg:riot_primal_precond}
    \begin{algorithmic}[1]                    

  \tiny   \STATE Start with $\mathbf{x}_{0}=\mathbf{x}_b$, $\mathbf{v}_0=\mathbf{0}$, $\mathbf{B}=\mathbf{L}_\mathbf{B}\mathbf{L}_\mathbf{B}^\top$, $\mathbf{P}'_0 = \mathbf{I}_n$, precond=true or false (logical), rotation=true or false (logical).
        \STATE Choose $m$ (number of samples for randomized SVD).
        \STATE Choose $k_f$ (number of outer loops).
        \FOR{$k=1, 2,\hdots, k_f$}         
         \STATE Integrate and store trajectory $H\left(\mathbf{x}_{k-1}\right)$.                 
          \STATE Compute and store $\mathbf{d}_{k-1}=H(\mathbf{x}_{k-1})-\mathbf{y}$.
           \STATE Let $\hat{\mathbf{A}}_0 \equiv \mathbf{L}_\mathbf{B}^\top\mathbf{H}_{k-1}^\top\mathbf{R}^{-1}\mathbf{H}_{k-1}\mathbf{L}_\mathbf{B} \in \mathbb{R}^{n \times n}$, where $\mathbf{H}_{k-1}$ is the tangent-linear and $\mathbf{H}_{k-1}^\top$ the adjoint model at $\mathbf{x}_{k-1}$.
    \vskip 2pt     
           \IF {\text{precond}}
             \IF {$k \ge 2$} 
\vskip 2pt
           \STATE $\mathbf{P}'_{k-1}= \mathbf{I}_n + \widetilde{\mathbf{V}}_{k-1} \left( \widetilde{\mathbf{\Lambda}}_{k-1}^{-\frac{1}{2}} - \mathbf{I} \right)\widetilde{\mathbf{V}}_{k-1}^\top$
\vskip 2pt
           \STATE  ${\mathbf{P}'_{k-1}}^{-1} = \mathbf{I}_n + \widetilde{\mathbf{V}}_{k-1} \left( \widetilde{\mathbf{\Lambda}}_{k-1}^{\frac{1}{2}} - \mathbf{I} \right)\widetilde{\mathbf{V}}_{k-1}^\top$
\vskip 2pt
                     \ENDIF
                        \STATE  $\mathbf{P}_k = \mathbf{P}'_0 \mathbf{P}'_1 \hdots \mathbf{P}'_{k-1}$
          \vskip 2pt     
           \STATE ${\mathbf{P}_k}^{-1} = {\mathbf{P}'_{k-1}}^{-1} \hdots {\mathbf{P}'_1}^{-1}{\mathbf{P}'_0}^{-1}$
                     \vskip 2pt     

           \STATE $\hat{\mathbf{A}}_k = \mathbf{P}_k^\top\mathbf{P}_k - \mathbf{I}_n + \mathbf{P}_k^\top \hat{\mathbf{A}}_0 \mathbf{P}_k$, 
       \ELSE
         \STATE $\hat{\mathbf{A}}_k = \hat{\mathbf{A}}_0$
                  \STATE $\mathbf{P}_k = \mathbf{P}'_0$
       \ENDIF
                             	   \STATE Draw $\mathbf{\Omega}_k\in\mathbb{R}^{n \times m}\sim\mathcal{N}(0,1)$.
	     \FORALL {$i\in\{1,2,\hdots,m\}$}
	    \STATE $\mathbf{\omega}_{i}=\mathbf{\Omega}_k(:,i)$
                  \IF {$\{k \ge 2 \text{ .AND. precond .AND. rotation\}}$} 
                   \STATE  Compute SVD $[ \mathbf{P}_1\widetilde{\mathbf{V}}_{1} ,\hdots,\mathbf{P}_{k-1}\widetilde{\mathbf{V}}_{k-1}]=  \mathbf{U}_k\mathbf{D}\mathbf{Z}_k^T$ (or use QR decomposition).
                    \STATE $\mathbf{\omega}_{i}={\mathbf{P}_k}^{-1}(\mathbf{I}_n - \mathbf{U}_k\mathbf{U}_k^T)\mathbf{\omega}_{i}$
                   \ENDIF
\ENDFOR
           \FORALLP {$i\in\{0,1,2,\hdots,m\}$}
                  \IF {$i = 0$}                     
                      \STATE Compute the gradient $\mathbf{b}\equiv\mathbf{P}_k^\top\mathbf{L}_\mathbf{B}^\top\mathbf{H}_{k-1}^\top\mathbf{R}^{-1}\mathbf{d}_{k-1}$.
                 \ELSE
                     \STATE $\mathbf{y}_{i} = \hat{\mathbf{A}}_k\mathbf{\omega}_{i}$               
                  \ENDIF
           \ENDFOR           

           \STATE Form $\mathbf{Y}=\left[\mathbf{y}_1,\mathbf{y}_2,\hdots,\mathbf{y}_{m}\right]$.
           \STATE Compute $\mathbf{Q}\in\mathbb{R}^{n \times m}$ orthonormal basis of $\mathbf{Y}$ from, e.g., $\text{QR}$ algorithm.
           \STATE Solve for $\mathbf{K}$ in $\mathbf{K}\mathbf{Q}^\top\mathbf{\Omega} = \mathbf{Q}^\top\mathbf{Y}$.
           \STATE Form eigendecomposition, $\mathbf{K}=\mathbf{Z}\mathbf{\Lambda}\mathbf{Z}^\top$ ($\mathbf{\Lambda} = diag(\lambda_i)$).
                      \STATE Compute  $\widetilde{\mathbf{V}}_{k}=\mathbf{Q}\mathbf{Z}$.
                      \STATE  $\widetilde{\mathbf{\Lambda}}_{k-1} = \mathbf{\Lambda} + \mathbf{I}$.
  \IF{ $\lambda_m>1$ }
         \STATE $\delta\mathbf{v}= -  \left (\sum_{i=1}^m \frac{1}{1+\lambda_i}\mathbf{u}_i\mathbf{u}_i^\top \right )\left ( \mathbf{b}+\mathbf{v}_{k-1} \right ) $
  \ELSE
                 \STATE $\delta\mathbf{v}= -\left ( \mathbf{I}-\sum_{i=1}^m \frac{\lambda_i}{1+\lambda_i}\mathbf{u}_i\mathbf{u}_i^\top\right )\left ( \mathbf{b}+\mathbf{v}^{k-1} \right ) $
  \ENDIF
          \STATE $\mathbf{v}_{k} = \mathbf{v}_{k-1}+\delta\mathbf{v}$
          \STATE $\mathbf{x}_{k} = \mathbf{x}_{k-1}+ \mathbf{L}_\mathbf{B}\mathbf{P}_k\delta\mathbf{v}$
        \ENDFOR 
    \end{algorithmic}
%
\end{algorithm}

\section{Numerical Experiments}
\label{num_exp}
\subsection{Lorenz-96 System}
\label{l96}
\subsubsection{Description}
The nonlinear dynamical Lorenz-96 system is frequently used to test new data assimilation methods, and is defined as follows \citep{Lorenz96}:

 \begin{align}
 \label{eq7}
\frac{d\mathbf{x}_i}{dt} = (\mathbf{x}_{i+1}-\mathbf{x}_{i-2})\mathbf{x}_{i-1}-\mathbf{x}_i+F,
\end{align}
where $\mathbf{x}_i$ represents the state of the system, $\mathbf{x}_{-1}=\mathbf{x}_{N-1}, \mathbf{x}_0=\mathbf{x}_N$, $\mathbf{x}_{N+1}=x_1$, and $F$ is a forcing term.

The dimension, $N$, and forcing, $F$, are free parameters of the Lorenz-96 system.  Often $F$ is set to $8$ and $N=40$ to ensure chaotic dynamics.  In our test $F=8$, but $N=300$ in order to increase the potential degrees of freedom (DOFs). 
 The initial true state, $\mathbf{x}_{0}$, was obtained from a spin-up simulation that was initialized with $\mathbf{x}_{-1}=[1, ...,1]$ and perturbation of $0.008$ at $\mathbf{x}_{N/2}$, then integrated forward in time for 144,000 time steps. 
In this experiment $\mathbf{x}_0$ was used as the truth in a series of synthetic observation inversions where the pseudo-observations are obtained by integrating the model from $\mathbf{x}_0$, i.e., $\mathbf{y}_\text{pseudo}=H(\mathbf{x}_0)$. In this case 100 observations are distributed randomly in space and time. 
Two assimilation windows are considered in our experiment, a 6-hour, 48-hour, 72-hour, and 96-hour windows with total integrations of $\Delta t = [0.05, 0.4, 0.6, 0.8]$. All four scenarios use a time-step of $\delta t = 0.01$.  The prior error variances are defined as $\sigma_b^2=(0.04F)^2 + (0.1 |\mathbf{x}_0 - F|)^2$ and the background error correlations as Gaussian functions with a length-scale of 1.5 grid cells. Observation error variances are defined as $\sigma_\text{obs}=(0.04F)^2 + (0.05|\mathbf{y}_\text{pseudo} - F|)^2$. Observation errors are assumed to be uncorrelated. The observations used in the assimilation are obtained by perturbing the pseudo-observations according to their error statistics: $\mathbf{y} = \mathbf{y}_\text{pseudo} + \mathbf{R}^{1/2} \mathbf{\epsilon}_o$, where $\epsilon_o\sim\mathcal{N}(0,1)$. Similarly, the prior used in the inversion is obtained by perturbing the truth: $\mathbf{x}_b = \mathbf{x}_0+ \mathbf{B}^{1/2}\epsilon_b$ where $\epsilon_b~\mathcal{N}(0,1)$.
\subsubsection{Results and Discussion}

\begin{figure*}[]
\centering
 \includegraphics[width=37pc]{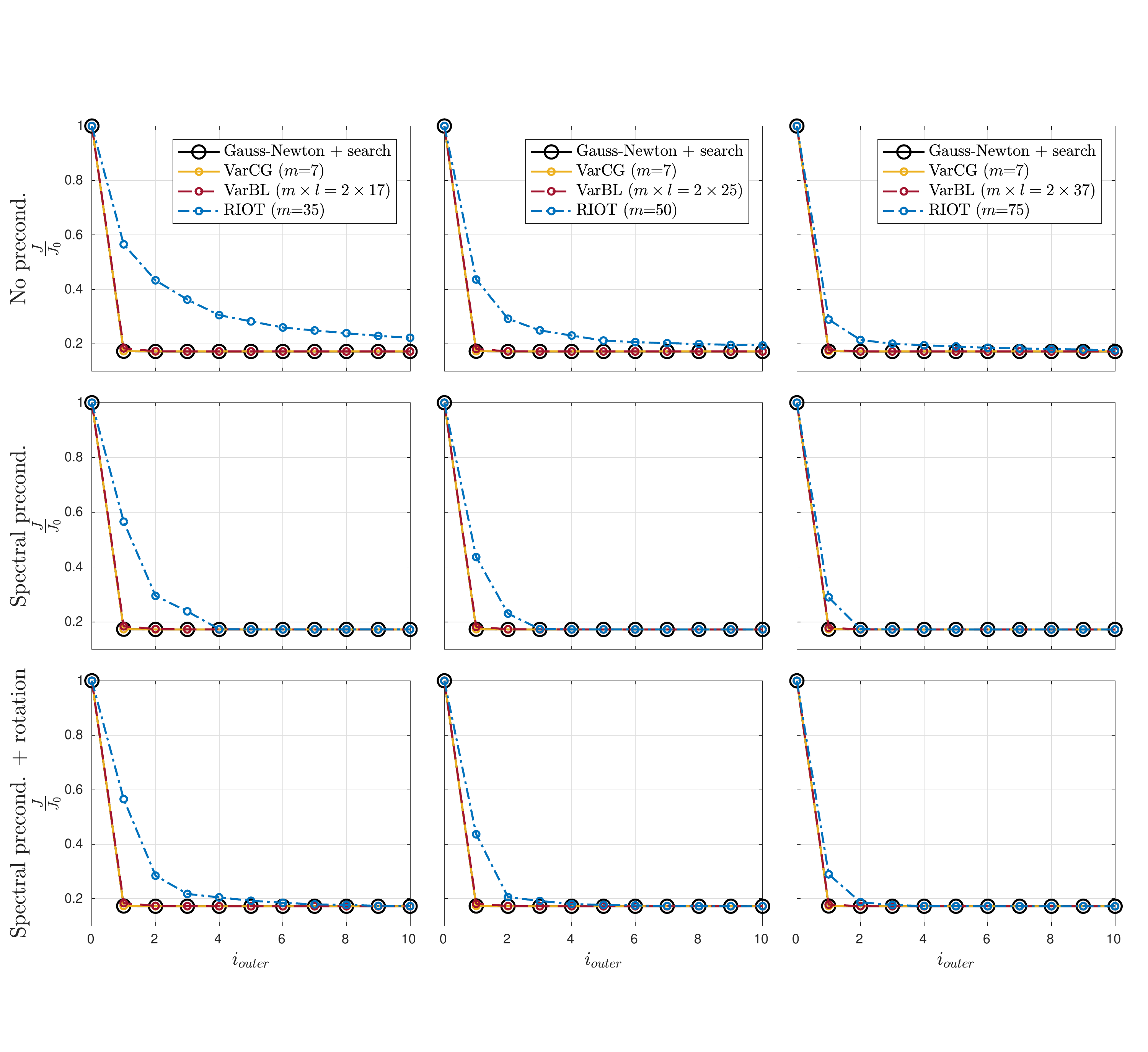}
\caption{Performance of the VarCG (yelow curve), VarBL (red curve) and RIOT (blue curve) mimimization algorithms for the 6-hour window L-96 problem. The figures show the non-quadratic cost function values (y-axis) at the beginning of each outer iteration (x-axis). Results are presented for different ranks of the Hessian approximation (from left to right column) and for different preconditioning approaches: no preconditioning (top row), spectral preconditioning (middle row), spectral preconditioning with rotation for VarBL and RIOT (bottom row). The parameter $m$ represents the number of samples in RIOT and the number of iterations in VarCG. For VarBL, $m$ is the number of iterations while $l$ is the number of samples (see algorithms \ref{alg:lanczos}, \ref{alg:Block_lanczos}, and \ref{alg:riot_primal_precond}). Results for a Gauss-Newton minimization (black curve) using an exact Hessian is also shown as reference.}
\label{fig:l96_6hrs_costf}
\end{figure*}

\begin{figure*}[]
\centering
 \includegraphics[width=37pc]{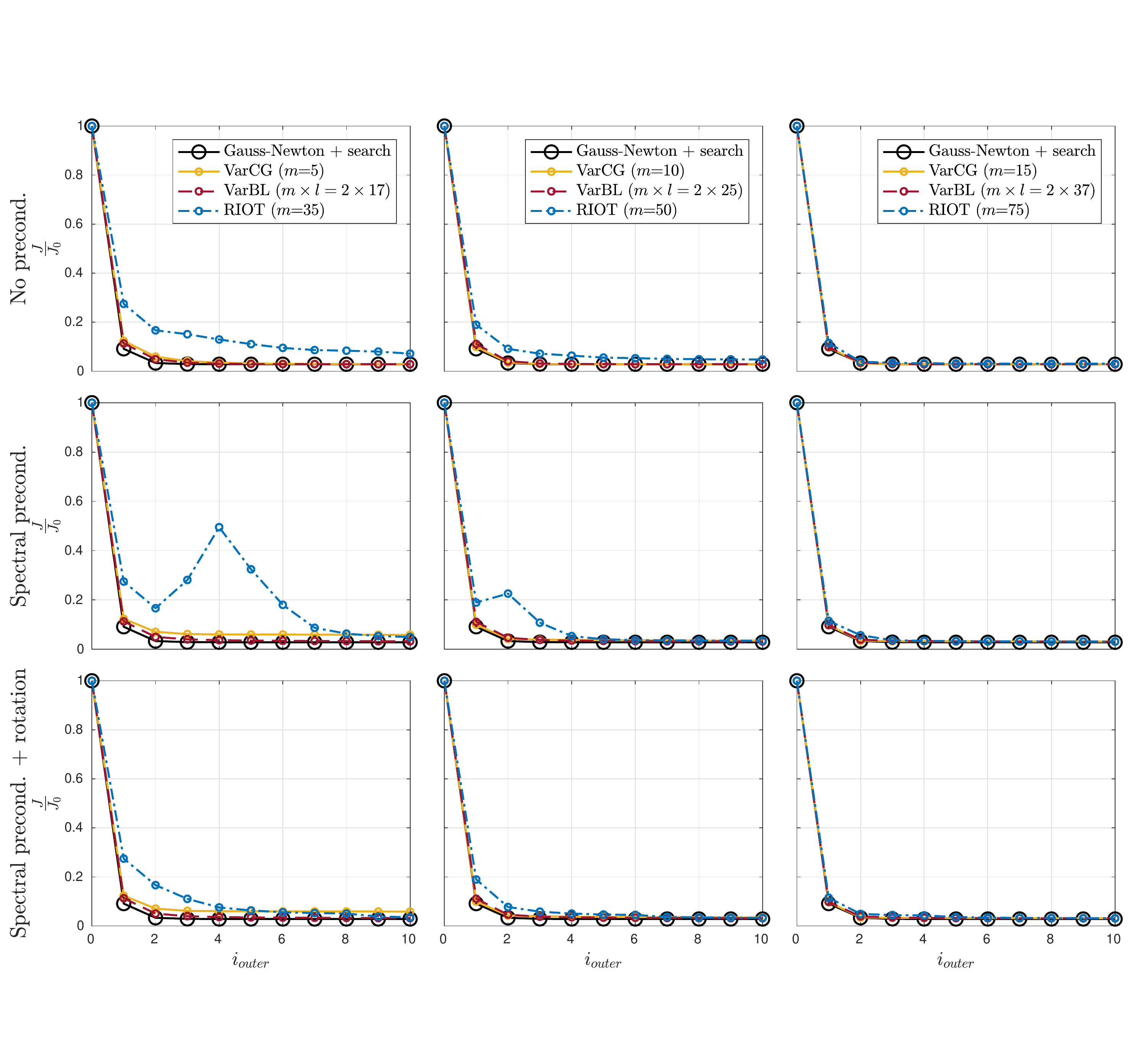}
\caption{Same as Fig. \ref{fig:l96_6hrs_costf} but for the 48-hour window L-96 problem.}
\label{fig:l96_48hrs_costf}
\end{figure*}

\begin{figure*}[]
\centering
 \includegraphics[width=37pc]{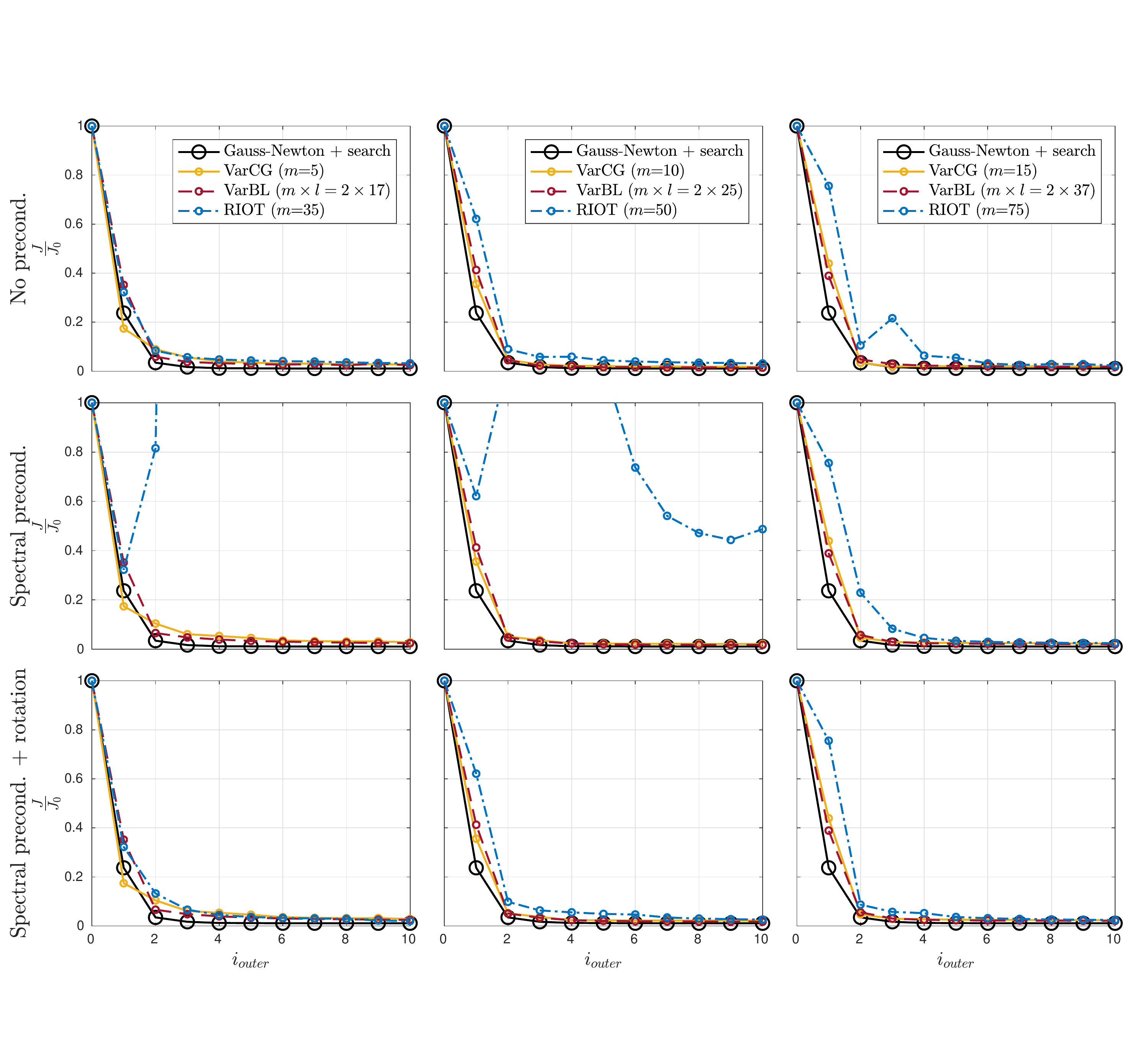}
\caption{Same as Fig. \ref{fig:l96_6hrs_costf} but for the 72-hour window L-96 problem.}
\label{fig:l96_72hrs_costf}
\end{figure*}

\begin{figure*}[]
\centering
 \includegraphics[width=37pc]{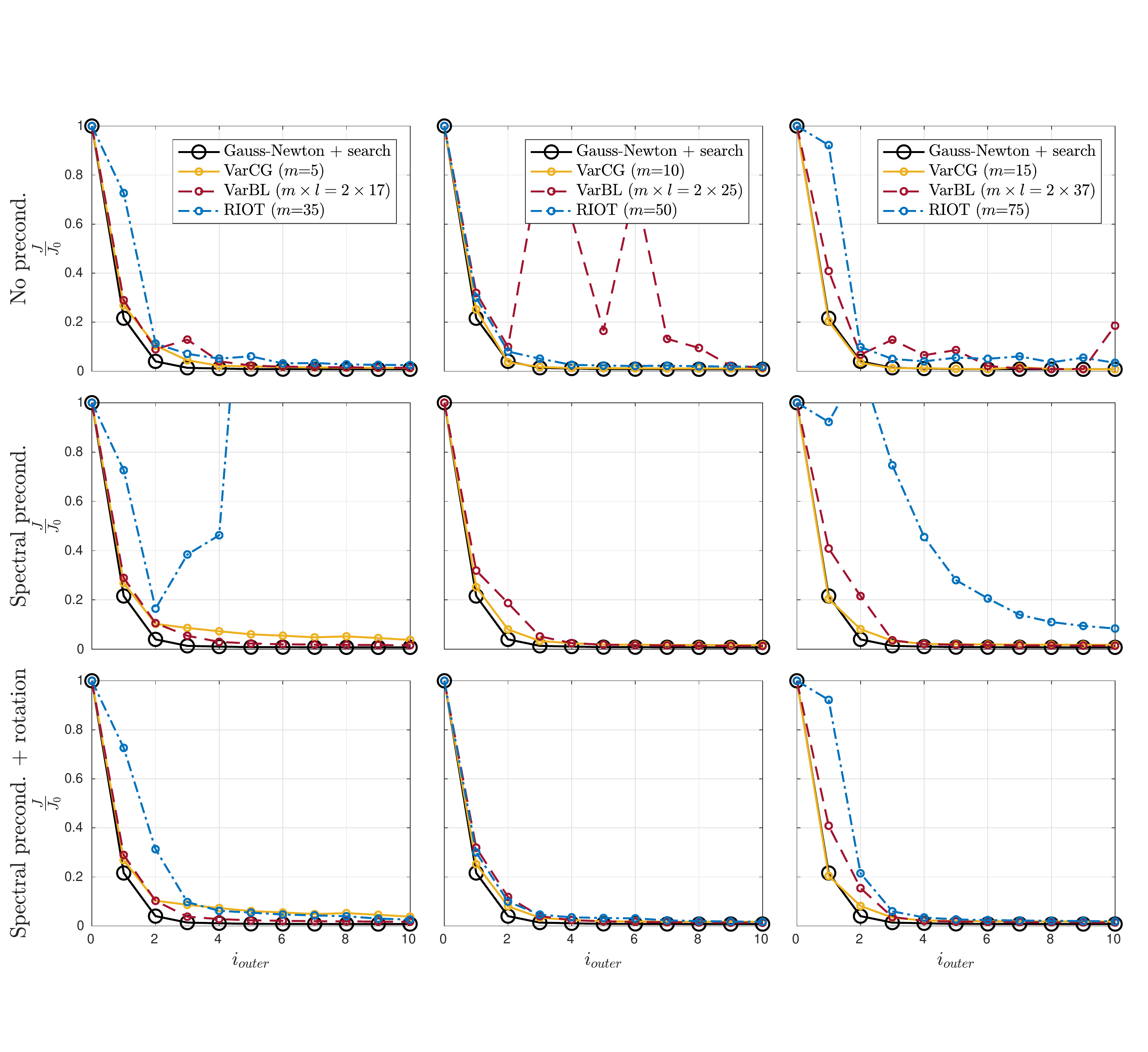}
\caption{Same as Fig. \ref{fig:l96_6hrs_costf} but for the 96-hour window L-96 problem.}
\label{fig:l96_96hrs_costf}
\end{figure*}

Figures \ref{fig:l96_6hrs_costf}, \ref{fig:l96_48hrs_costf}, \ref{fig:l96_72hrs_costf} and \ref{fig:l96_96hrs_costf} show the non-quadratic cost function values at the beginning of each outer iteration for VarCG, VarBL, and for RIOT, for the 6-hour, 48-hour, 72-hour, and 96-hour L-96 problem, respectively.  Also included is a reference solution computed using Gauss-Newton and the exact Hessian.  For the 6-hour experiment, the VarCG minimization converges in the first outer iteration (i.e., after only 7 inner iterations ), which reflects the strongly linear behavior of the system. The performance of the RIOT algorithm is sensitive to both the number of samples used in the RSVD ($m$) and the preconditioning. For example, the top row of Fig.~\ref{fig:l96_6hrs_costf} shows that when using 35 samples in the RSVD the minimum of the cost function is not reached after 10 outer iterations, whereas it is reached when 75 samples are used.

Using the spectral preconditioning in RIOT significantly improves the convergence rate of the minimization in the 6-hour case, especially when using 35 or 50 samples in the RSVD. For instance, the spectral preconditioning technique enables RIOT with 35 samples to convergence in about 6 outer iterations, similar to the non-preconditioned RIOT case with 75 samples.  Preconditioning RIOT thus mitigates the need to increase the ensemble size in the RSVD. For the 6-hour window problem, the rotation approach does not improve the results compared to those obtained with spectral preconditioning alone. 

The same is not true for the more non-linear cases (48-hour to 96-hour windows), where the spectral preconditioner produces spurious increases in the cost function, in some cases resulting in the divergence of the RIOT algorithm (72-hour and 96-hour cases). As explained in section \ref{methods}, the preconditioner reduces the condition number of the Hessian, which can degrade the performance of the RSVD in approximating the dominant eigenvectors by increasing the interference of weaker modes. Introducing the rotation step allows one to efficiently filter out part of that noise by reorienting the RSVD input samples in the direction of the dominant Hessian modes. The preconditioning and rotation approaches appear to be beneficial to VarBL as well in strongly non-linear regimes (see 96-hour window results). 

In all L-96 experiments, combining the spectral preconditioner with the rotation method in the RSVD enables RIOT to converge to the Gauss-Newton reference solution either notably more rapidly (for the weakly non-linear case) or as rapidly (highly non-linear case) as the non-preconditioned case.  Conversely, the spectral preconditioner applied to VarCG does not improve the performance of the minimization. As noted by \citet{Tshimanga2008}, using a spectral preconditioner that implicitly assumes exact eigenvectors computation in place of the more consistent Ritz preconditioner can result in poorer performances than obtained from a non-preconditioned VarCG minimization. The impact of using the Ritz preconditioner instead of the spectral preconditioner is beyond the scope of this paper. Note that tests conducted with the rotation modification only in RIOT (i.e., without spectral preconditioning) did not show any improvement of the convergence of the minimization except for the almost linear case corresponding to the 6-hour window (results are not shown here).

In conclusion, although for a given rank of the Hessian approximation (i.e., number of samples in VarBL and RIOT, or number of iterations in VarCG) VarCG provides better minimization performances than either RIOT or VarBL, the preconditioning techniques applied to the latter methods ensure  convergence  in both linear and non-linear regimes, even using a moderate number of samples for the RSVD Hessian approximations. A noteworthy characteristic of the spectral preconditioning with rotation used in RIOT or VarBL is its stability compared to all other approaches. This much desirable property of our preconditioning technique is clearly observed in the strongly non-linear cases (Fig. \ref{fig:l96_72hrs_costf}\ and \ref{fig:l96_96hrs_costf}), where both RIOT and VarBL exhibit erratic behavior of the minimization with sudden increases in the cost function value when no preconditioning is applied. Lastly, as discussed in section \ref{methods}, both VarBL and RIOT enable highly parallel implementations of incremental 4D-Var, the computational wall-time of the minimization has the potential to be drastically reduced compared to the sequential VarCG method, provided enough nodes are available to run many instance of the TL and AD models in parallel. In the next section we illustrate this aspect of the performance of RIOT in a more practical context using two large-scale atmospheric source inversion problems.


\subsection{Application to an Atmospheric Source Inversion Problem}
\label{smallBC} In this section we further analyze the properties of the RIOT algorithm in a more realistic framework by considering two large-scale atmospheric source inversion problems. Both experiments consist of optimizing regional black carbon emissions using aircraft measurements from the ARCTAS-CARB (Arctic Research of the Composition of the Troposphere from Aircraft and Satellites in conjunction with the California Air Resources Board) campaign.  

\subsubsection{Performance Evaluation for a low DOFs problem}
\label{smallBC}
\begin{figure*}[]
\centering 
 \includegraphics[width=20pc]{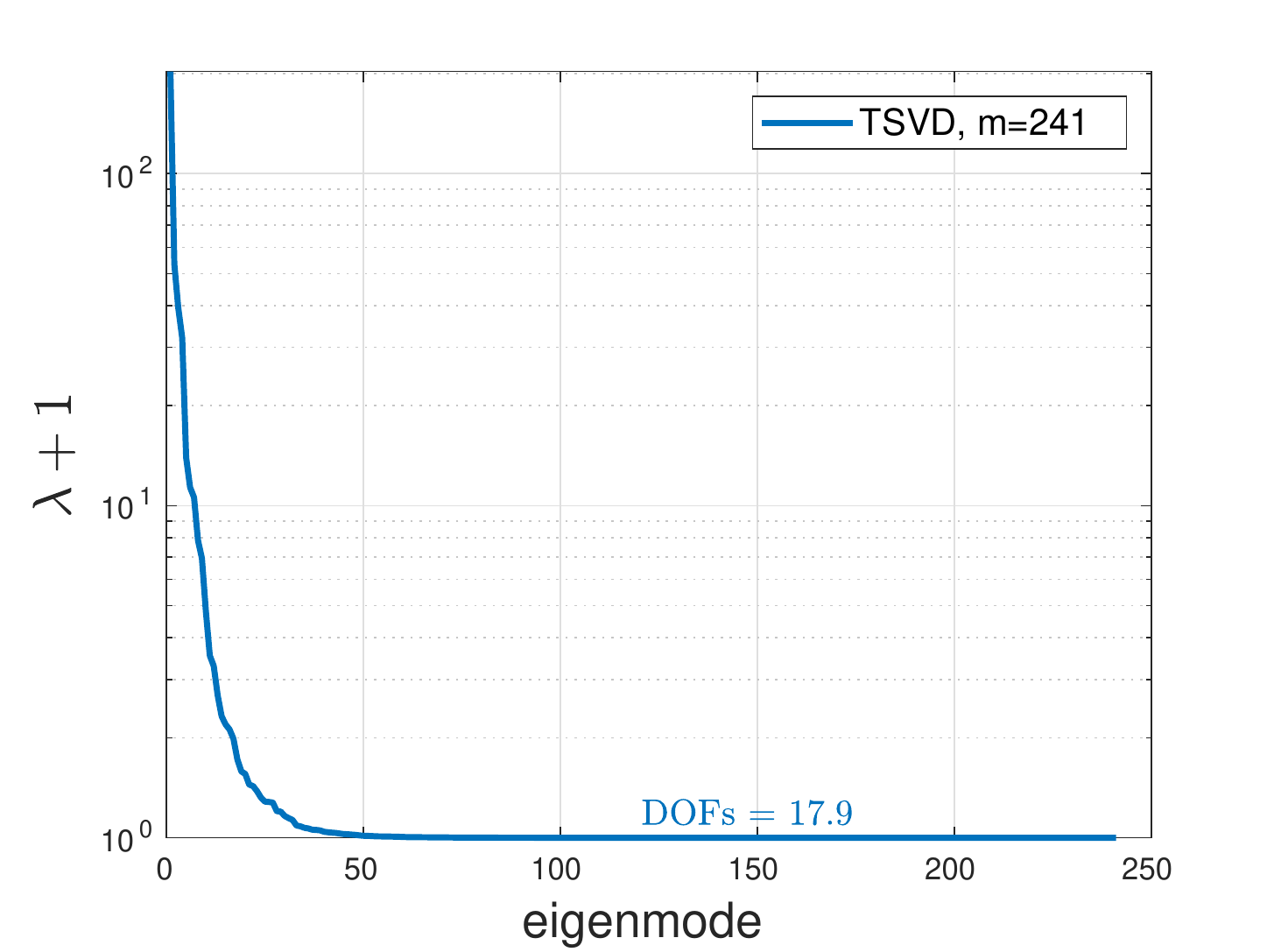}
\caption{Eigenspectrum of the prior-preconditioned Hessian for the first outer iteration of the 4D-Var minimization for the small BC problem. The corresponding DOFs ($\mathrm{DOFs}=\sum_i\lambda_i/(1+\lambda_i)$) is indicated.}
\label{fig:bcsmall_eigs}
\end{figure*}

\paragraph{Description}
The small black carbon source inversion problem (referred to as small BC thereafter) uses real data from an SP2 device measuring black carbon on a single flight from ARCTAS-CARB on 22 June, 2008 [ref]. 
The experiment is identical to FINN$_\text{STD}$, as described in \citet{Guerrette17}.  The prior inventories are derived from FINN biomass burning (BB) emissions and NEI 2005 anthropogenic (ANT) emissions.  Scaling factors for the BC emissions resolved hourly and across grid cells are optimized assuming a prior value of 1. A factor of 3.8 uncertainty for BB and factor of 2 for ANT are prescribed, both with a 36 km horizontal correlation length scale. The results here are slightly different from those described in \citet{Guerrette17}, because since that paper was published errors in the trajectories associated with the solar radiation treatment have been diagnosed and fixed.  Improved symmetry between TL and AD used for this study give more accurate results than those in \citet{Guerrette17}. The previously existing asymetry was severe enough to prevent RSVD from converging effectively, but was less detrimental to VarCG. The total number of observations considered for the small BC problem is 241.  While the total dimension of the control vector was 299568, the number of grid columns with active emission values was 121912 and the number of grid cells with significant emission values is far fewer. As a result, the small BC problem has a low DOFs of 17.9 (see Fig. \ref{fig:bcsmall_eigs}). For a problem of those dimensions, the exact computation of both the posterior mean emissions and the posterior error covariance matrices is possible. Therefore the approximations obtained from the different minimization algorithms can be evaluated against the exact solutions. For the small BC problem, the VarBL algorithm was not tested, therefore only results obtained with the VarCG and RIOT methods are presented and discussed. As described in \citet{Guerrette17}, a lognormal transformation is applied to the emission control vector in order to enforce positiveness of the solution. Therefore, although the transport model is a linear solver, non-linearities are introduced in the resulting forward model that requires the use of non-quadratic minimization algorithm such as incremental 4D-Var.
\paragraph{Results and Discussion}
  In this section we first analyze the behavior of VarCG and RIOT for the inner loop minimization. In particular, we assess the impact of the adaptive approach used in RIOT in order to guarantee statistical optimality of the approximated posterior increment and posterior covariance matrices. 
  
\begin{figure*}[]
\centering
 \includegraphics[width=15pc]{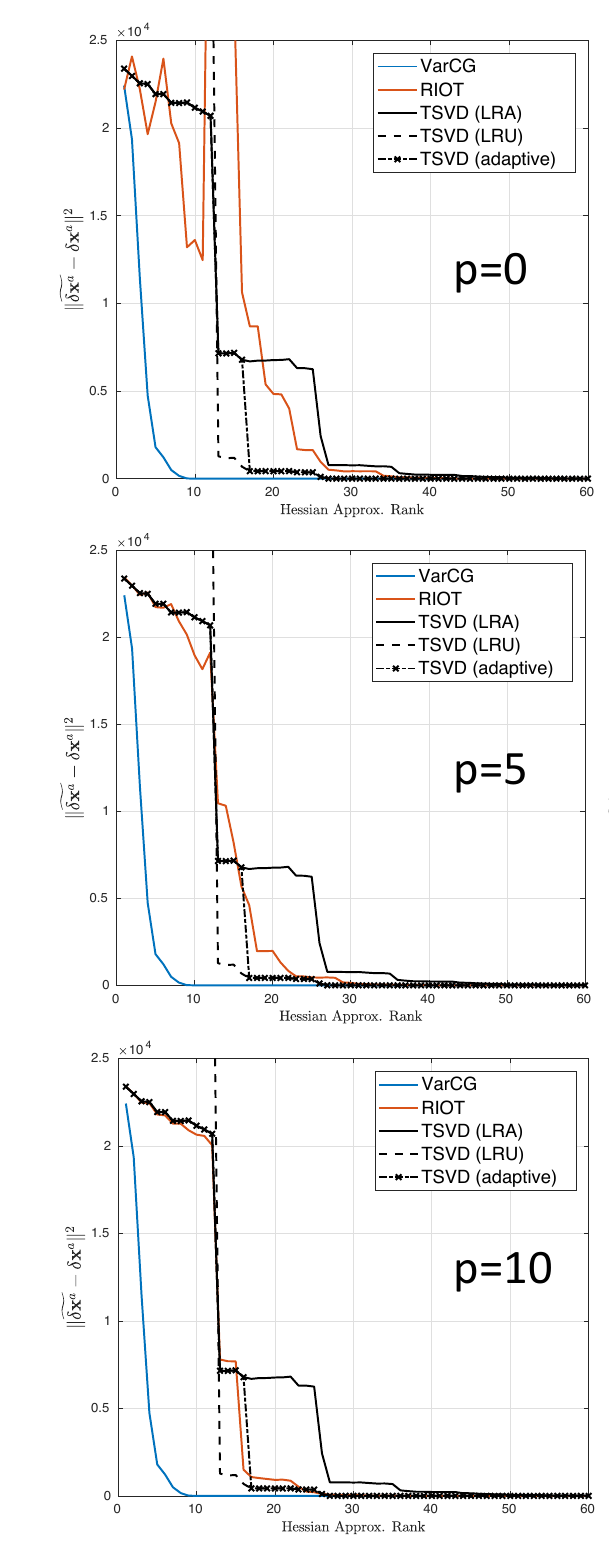}
\caption{Increment approximation error (Euclidian norm) (y-axis) for the first inner loop of the small BC problem as a function of the rank of the Hessian approximation (x-axis) for different methods: VarCG (blue line), RIOT (red line), TSVD LRA (black solid line), TSVD LRU (black dashed line) and TSVD adaptive (black doted line) updates. From top to bottom, the RIOT algorithm uses an oversampling parameter in the randomized SVD with values p=0, p=5 and p=10, respectively.}
\label{fig:bcsmall_1rstouter}
\end{figure*}

 Figure \ref{fig:bcsmall_1rstouter} shows the squared Euclidian norm of the error in the posterior increment as a function of the rank of the Hessian approximation for the first outer iteration of the 4D-Var minimization. The rank of the Hessian approximation $k$ is equal to the number of iterations in the case of VarCG, and to the number of approximated eigenvectors retained in the RSVD Hessian approximation in the case of RIOT.  Figure \ref{fig:bcsmall_1rstouter} also presents the performance of RIOT for different values of the oversampling parameter $p$, whose value defines the total number of samples $m = k+p$ used in the RSVD approximation. 

 Results from a Hessian approximation with exact truncated eigendecomposition (TSVD), where only the first $k$ eigenvectors are retained are also shown for each rank $k$.
The RIOT increments show very noisy features in the case where no oversampling is used in the  RSVD ($p=0$). It appears that a small overampling parameter of $5$ is enough to obtain an increment similar to the adaptive approach applied to the exact truncated SVD approximation, while $p=10$ yields almost perfect agreement. This is consistent with \citet{Halko11}, in which oversampling parameters of $5-10$ were sufficient to ensure good convergence properties in most cases. As expected from their statistical properties (see \citet{bous18}), TSVD-LRA is associated with smaller increment errors than TSVD-LRU from rank 1 to 12, which correspond to prior-preconditioned Hessian eigenvalues $\lambda_i>1$ (not shown), while TSVD-LRU outperforms TSVD-LRA for higher ranks. One observes that the adaptive approach applied to TSVD provides a solution that is close to the optimal (minimum increment error) one. Nevertheless, the adaptive selection of the LRU update occurred 3 ranks too late (at rank 15) compared to the optimal transition point. This can be explained by both the statistical nature of the criteria, that is, the optimality is valid in average, not necessarily for a particular realization of the prior and observations, and the fact that the optimal transition depends (even in average) on the full spectra of the Hessian. More specifically, while in average TSVD-LRU will always correspond to smaller errors than TSVD-LRA for truncation indices associated with eigenvalues $\lambda_i<1$, the reverse does not necessarily hold for $\lambda_i>1$. As a result, in general the rank associated with the optimal transition from TSVD-LRA to TSVD-LRU can be smaller than the one inferred from the RIOT adaptive criteria. Those results suggest another, deterministic method for the adaptive approach. Assuming a maximum of $m_\text{max}$ samples can be afforded for each inner iteration in RIOT, once the rank-$m_\text{max}$ RSVD approximation of the Hessian is available it is possible to perform two parallel integrations of the non-linear forward model, one using the TSVD-LRA and the other the TSVD-LRU increment update, in order to evaluate the non-quadratic cost function J. The increment update corresponding to the smallest J would then be selected, and the same procedure applied to the next inner loop. In practice this would require an additional (possibly parallel) integration of the non-linear model compared to the current RIOT algorithm. 

For a given rank $k$, the VarCG method significantly outperforms both RIOT and the TSVD approach, reaching convergence after only 10 inner iterations. This can be explained by the fact that the CG algorithm ensures minimum error residuals at each iteration \citep{Golub89}. In particular, a key difference between the RSVD-based or TSVD Hessian approximations and the implicit Hessian approximation obtained from the CG minimization is that the latter method exploits the gradient information (i.e., $\mathbf{b}$ in the right-hand side of (\ref{eq4}).) While RIOT utilizes approximations of the  modes of the Hessian in order of decreasing eigenvalue, CG tends to approximate Hessian modes that provide significant contribution to the reduction of the error residual. In this experiment, since only about $25$ samples were needed for RIOT to converge (the convergence criteria being defined as $|| \widetilde{\delta \mathbf{x}^a}-\delta \mathbf{x}^a ||^2 < 10$), all randomized SVD integrations were run in parallel. This yielded about an order of magnitude decrease in wall-time compared to the 10 iterations required for the VarCG algorithm. 

\begin{figure*}[]
\centering 
 \includegraphics[width=40pc]{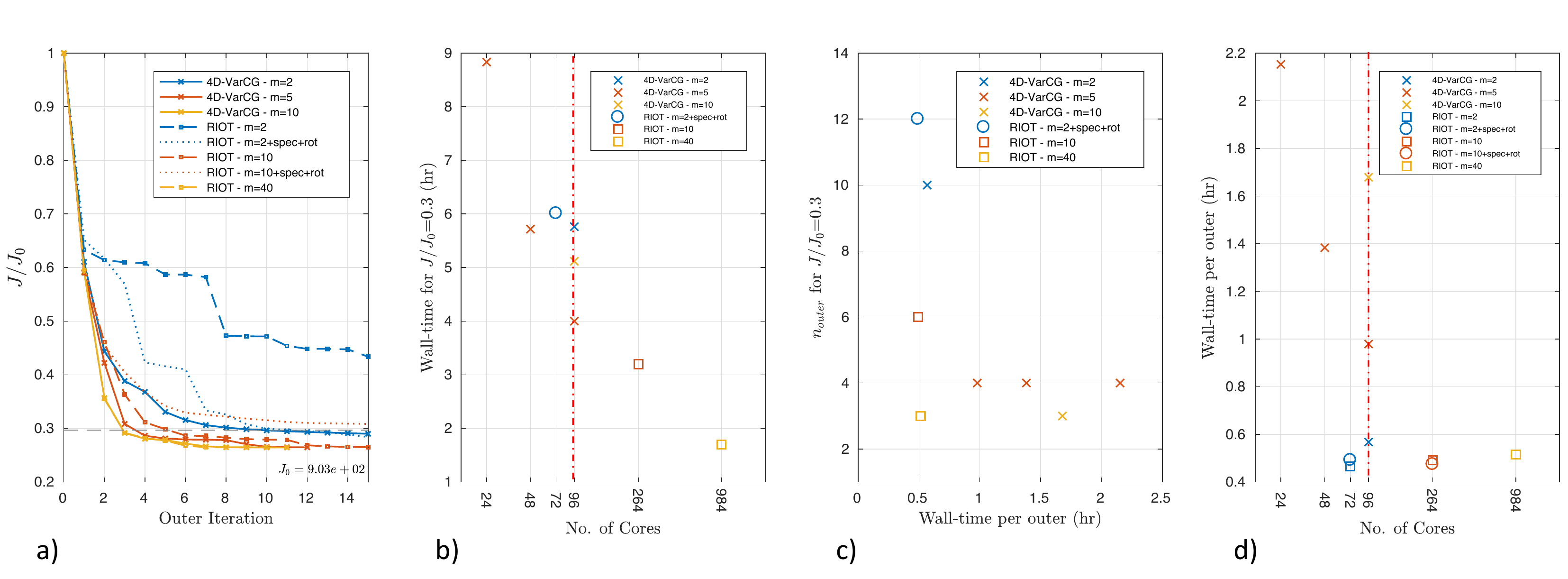}
\caption{Convergence performance for different configurations of RIOT and VarCG for the small BC problem. The VarCG experiments differ by the number of inner-loop iterations (m) and/or the number of cores used per model integration. The RIOT experiments differ by the rank of the Hessian approximation (m) and by the use of the spectral (spec) preconditioning and/or the rotation method (rot). a) Relative cost function value (fraction of initial value) (y-axis) as a function of the outer iteration number (x-axis); b) wall-time to reach a cost function reduction of 70\% (y-axis) and corresponding number of cores per model integration (x-axis); c) number of outer iterations to reach a cost function reduction of 70\% (y-axis) and corresponding wall-time (x-axis); d) wall-time per outer iteration (y-axis) and corresponding number of cores per model integration (x-axis). The vertical red dotted line shows the maximum number of cores above which the MPI parallelization of the WRF-DA system does not reduce the computational wall-time of the integration. }
\label{fig:bcsmall}
\end{figure*}

 Figure \ref{fig:bcsmall} summarizes the performance of different VarCG and RIOT configurations in terms of cost function reduction, number of cores used per model integration and wall-time. Two metrics are considered, a local one based on computational performances per outer iteration, and a global one associated with performances for the entire minimization. A total of 15 outer iterations was used in those experiments, after which a minimization is considered to have converged if its final cost function value is within $3$\% of the minimum cost function across all methods. The selected threshold corresponds to the upper bound of cost functions values among all configurations that plateaued after 15 iterations. In terms of number of outer iterations required to reach convergence, RIOT with 40 RSVD samples ($m=40$) yields similar performance as the VarCG algorithm using 10 inner iterations ($m=10$). While VarCG with $m=10$ uses 96 cores per model integration, RIOT with $m=40$ uses about 10 times as many cores due to the 40 parallel model integrations running on 24 cores each (Fig. \ref{fig:bcsmall} d)). However, RIOT with $m=40$ enables a threefold decrease in wall-time per outer iteration compared to VarCG with $m=10$ (Fig. \ref{fig:bcsmall} d)). It should be noted that the wall-time performance of the MPI parallelization is intrinsically limited by the increase in communications relative to useful computations as more cores are used. In the case of the WRF model integrations, the robustness of those communications breaks down for more than 96 cores on the small regional model domain considered. Therefore, RIOT provides an additional level of parallelism that allows one to increase the wall-time performance beyond that obtained from a model-level MPI parallelization (Fig. \ref{fig:bcsmall} b)). 
 
 The effect of preconditioning in RIOT is also shown in Fig. \ref{fig:bcsmall}. A significant increase in the rate of convergence of the minimization is obtained when using the spectral preconditioner with rotation (see Sec. \ref{methods}) for RIOT with $m=2$, while the preconditioning technique slightly degrades the performance of the minimization compared to no preconditioning when $m=10$. For $m=40$ (not shown), the degradation obtained from preconditioning RIOT was even larger. These results illustrate two competing effects of the spectral preconditioning: reduction of the condition number that degrades the performance of RSVD methods and filtering of the Hessian modes resolved in previous outer iterations that enables the RSVD to focus on the remaining eigenmodes. For $m=2$, the rate of decrease of the preconditioned Hessian eigenspectra remains fast enough across the first outer iterations to obtain an accurate RSVD approximation of the unresolved modes. As the rank of the Hessian approximation is increased the eigenspectrum of the preconditioned Hessian tends to be flatter, which can generate noisy RSVD approximations. The degradation in the accuracy of the RSVD Hessian estimate can then offset the benefits obtained from filtering out the previously resolved modes using the spectral preconditioning. Therefore, as a general rule, the use of a spectral preconditioner in RIOT will tend to be beneficial for a Hessian with fast decreasing eigenspectrum and when the rank of the Hessian approximation is small (compared to the number of eigenvalues with significant gaps). Figure \ref{fig:bcsmall_eigs}, which shows the eigenspectrum of the prior-preconditioned Hessian in the first outer iteration of the minimization, further supports this interpretation. Indeed, when the rank of the Hessian approximation is very small (e.g., $m=2$) compared to the maximum eigenvalue indice associated with fast decreasing spectrum ($\approx 40$), one obtains a clear benefit from using the preconditioner in RIOT. Conversely, when the rank of the Hessian approximation corresponds to eigenvalue indices located in the flat regime of the spectrum (e.g., $m=40$), a non-preconditioned RIOT will outperform a preconditioned one.
  
\begin{figure*}[]
\centering 
 \includegraphics[width=40pc]{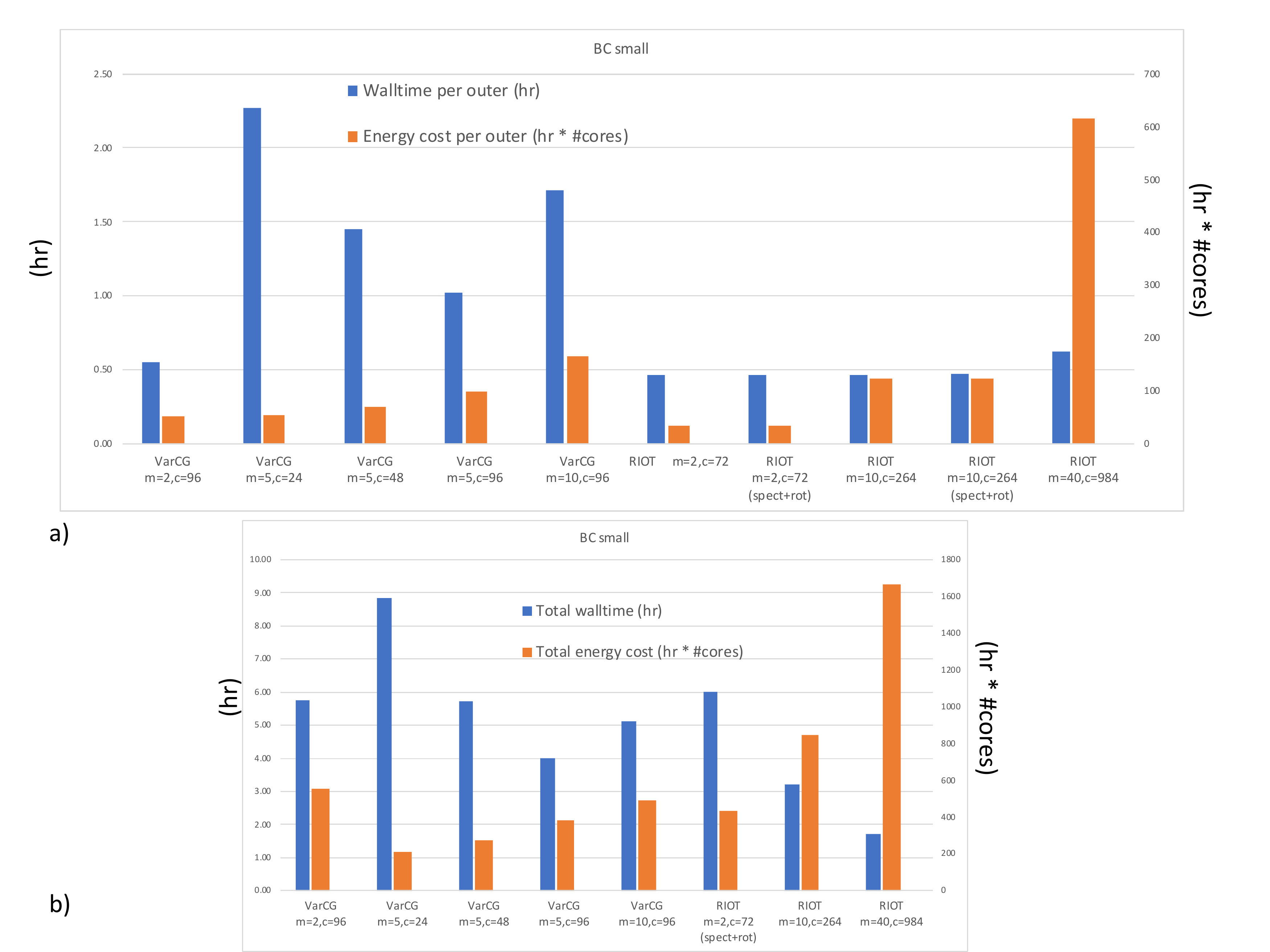}
\caption{Wall-time and energy cost performance for different configurations of RIOT and VarCG for the small BC problem. The number of inner iterations for VarCG or the rank of the Hessian approximation for RIOT are denoted by $m$, while the number of cores used per model integration is denoted by $c$. a) Wall-time (hours) and energy cost (hours$\times$cores) per outer iteration; b) wall-time (hours) and energy cost (hours$\times$cores) to reach a cost function reduction of 70\%}
\label{fig:bcsmall_chart}
\end{figure*}

 Figure \ref{fig:bcsmall_chart} summarizes the performances of the different minimization configurations in terms of both the wall-time and a metric directly related to energy consumption. The energy-related metric chosen here is the product of the number of cores simultaneously used and the wall-time, either per outer iteration (Fig. \ref{fig:bcsmall_chart} a)) or per converged minimization (Fig. \ref{fig:bcsmall_chart} b)). In term of per outer iteration cost, the RIOT configurations with 2 samples ($m=2$) performed best in term of both walltime and energy cost, followed by VarCG with $m=2$ (i.e., two inner-loop iterations) and $c=96$ cores per model integration. However, when considering the whole minimization (Fig. \ref{fig:bcsmall_chart} b)), VarCG configurations with more inner iterations ($m=5,\,10$) and at least 48 cores performed better than VarCG with $m=2$ iterations. Although it required more than 4 times as much energy demand as the fastest VarCG configuration ($m=5,\,c=96$), RIOT with $m=40$ samples allowed to reduce the walltime by more than a factor of 2 compared to the latter, making this algorithm an attractive approach when energy cost is not a significant limitation of the system considered.

In addition to the performance of the 4D-Var minimization, it is useful to compare the computational efficiency of RIOT and VarCG for estimating the posterior error covariance matrix, or, equivalently, the inverse Hessian of the cost function at its minimum. As described in Sec. \ref{methods}, another useful application of Hessian approximations is the preconditioning of the inner loop minimizations in incremental 4D-Var.

\begin{figure*}[]
\centering
 \includegraphics[width=25pc]{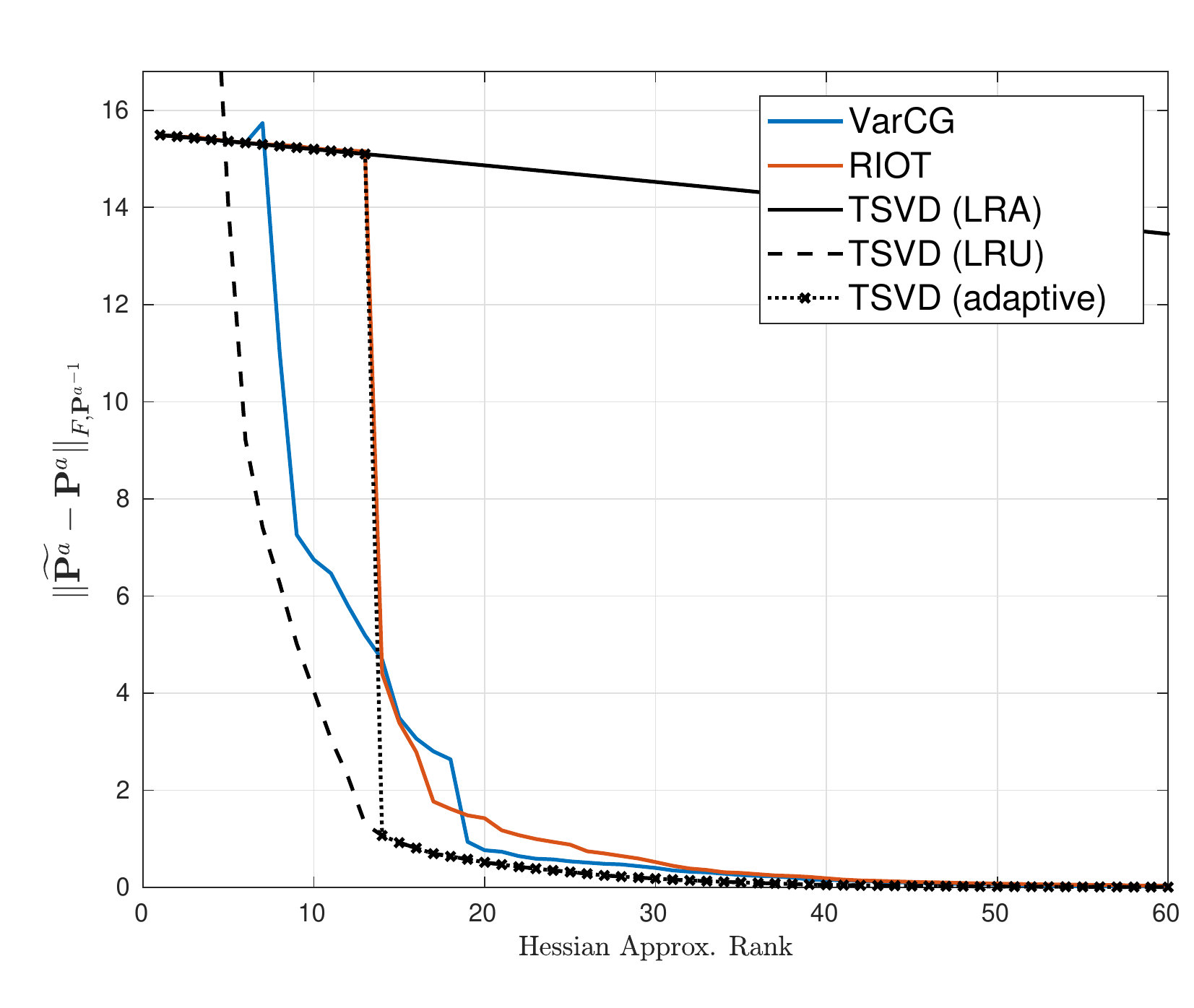}
\caption{Posterior error covariance approximation error (Pa-normalized Frobenius norm) as a function of the rank of the Hessian approximation for different methods: VarCG (blue line), RIOT (red line), TSVD LRA (black solid line), TSVD LRU (black dashed line), and TSVD adaptive (black dot line) approaches. An oversampling parameter of 5 ($p=5$) is used in the randomized SVD approximation. }
\label{fig:bcsmall_error}
\end{figure*}

Figure \ref{fig:bcsmall_error} shows the Frobenius norm of the relative error approximation in the posterior error covariance matrix \citep{bous18} as a function of the rank of the Hessian approximation for different methods. The rank of the approximation refers to either the number of iterations in the case of CG, the number of RSVD samples in the case of RIOT, or the number of eigenmodes retained in the case of TSVD. The posterior covariance approximation error for TSVD is shown for both the LRA and LRU updates, as well as for the adaptive method that is used in the RIOT algorithm (see Sec. \ref{methods}). The TSVD results clearly illustrate the importance of using the adaptive approach to define the posterior covariance approximation, with a fast decreasing error for LRU reaching a magnitude about 15 times smaller than that of the LRA estimates for ranks greater than 13 (where the adaptive transition to the LRU approximation occurs). However, it appears that the LRA estimate yields larger approximation errors than LRU even before the rank-13 approximation (that is, for ranks greater than 5). As explained in  \citet{bous18}, the superiority of the truncated LRU approximation over LRA is guaranteed only when the last mode used in the approximation is associated with an eigenvalue smaller than one. In the case where the smallest eigenvalue in the approximation is greater than one, using the LRA approximation may lead to a suboptimal choice, as we find here. The efficiency of the adaptive method, that is, its ability to select the optimal rank for the transition to the LRU approximation depends on the characteristic of the Hessian spectra itself. In particular, the adaptive approximation is more likely to generate more optimal posterior covariance estimates for fast decreasing eigenspectra than for relatively flat eigenspectra.
In practice, for a given inverse problem (e.g., data assimilation for NWP, atmospheric source inversion), one might adapt the LRA to LRU transition criteria to maximize the performance of the approximations based on the characteristics of the Hessian eigenspectra.

Figure \ref{fig:bcsmall_error} shows that both the RIOT and VarCG methods provide approximations to the posterior error covariance that are close to optimality for ranks greater than 20. We note that the VarCG approximation is associated with smaller errors than TSVD adaptive or RIOT from rank 8 to 12, which is due to an earlier transition to the LRU approximation in the case of VarCG (not shown here). This improvement is therefore an artifact of the underestimation of the true eigenvalues for indices 8 to 12 by the Lanczos-CG algorithm. The fact that the Lanczos method tends to approximate the extrema of the eigenspectrum (i.e., maximum and minimum eigenvalues) is well-known  \citep{Golub89}. The underestimation of the smallest approximated eigenvalues for a given rank $k<<n$ that we obtain here is thus consistent with the theory. For indices greater than 12, the exact eigenvalues quickly reach a minimum plateau (not shown here), which greatly mitigates the above-mentioned effect and explains the fast convergence of the VarCG posterior error covariance approximations to the TSVD approximations.

Finally, it is worthwhile to note that for all methods convergence to the exact posterior error covariance matrix is obtained for a rank 40 only. While it requires 40 iterations in the case of VarCG, it requires 40 parallel TL and AD integrations in the case of RIOT (i.e., a wall-time equivalent to one VarCG iteration). Thus, provided enough cores are available to perform all parallel RSVD integrations at once, RIOT allows one to reduce the wall-time cost of the approximation by a factor 40. This should be contrasted with the posterior increment approximation, for which the VarCG algorithm is more competitive due to its ability to optimally exploit gradient information and to enforce minimum residual error at each iteration through the conjugacy property. 
  
\subsubsection{Scalability for Problems With Larger DOFs}
\begin{figure*}[]
\centering 
 \includegraphics[width=20pc]{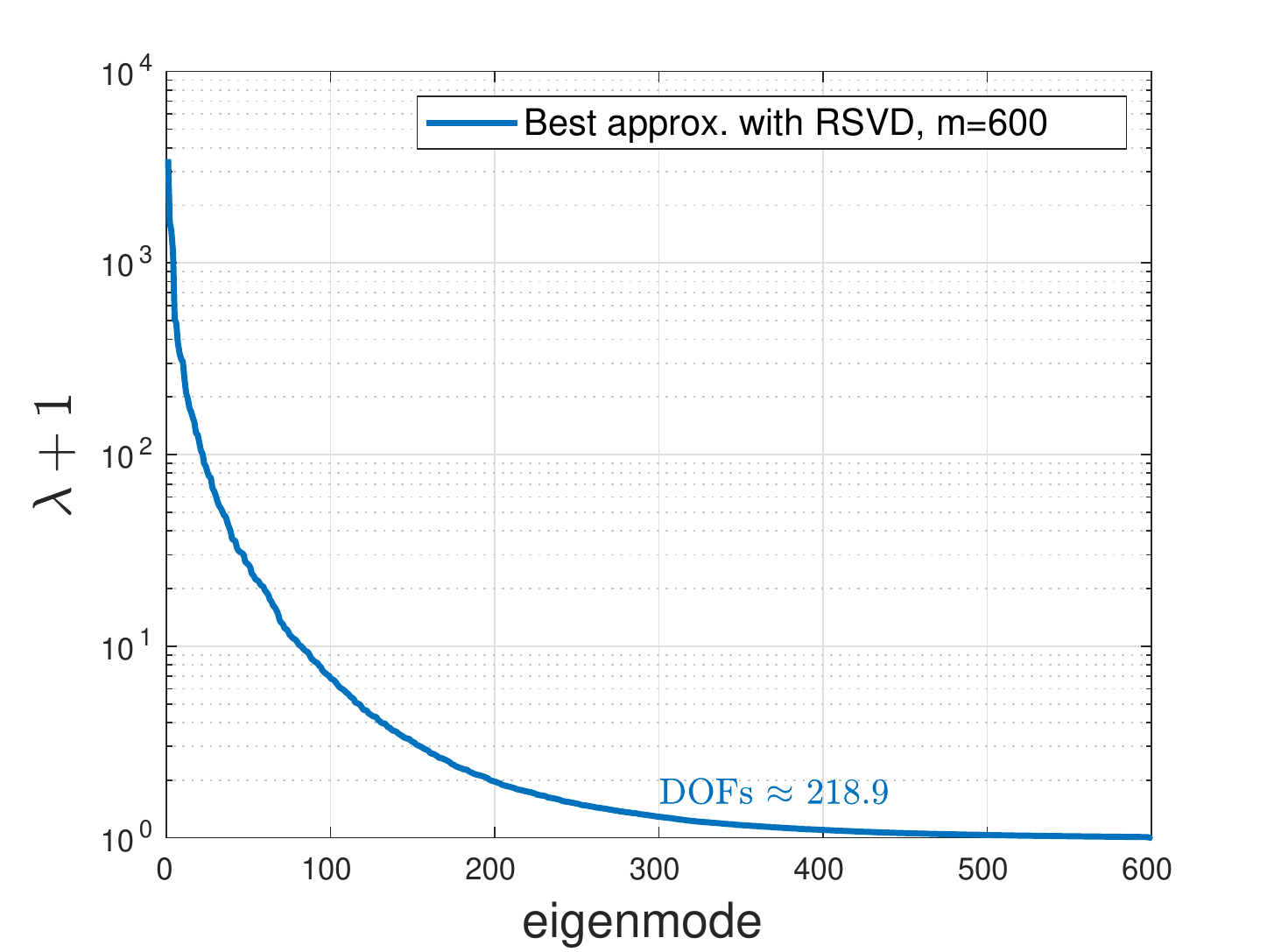}
\caption{Eigenspectrum of the prior-preconditioned Hessian for the first outer iteration of the 4D-Var minimization for the large BC problem. The corresponding DOFs ($\mathrm{DOFs}=\sum_i\lambda_i/(1+\lambda_i)$) approximated using 600 RSVD samples is indicated.}
\label{fig:bclarge_eigs}
\end{figure*}

\paragraph{Description}
The large black carbon source inversion problem (referred to as large BC thereafter) is designed to test the scalability of the RIOT algorithm for problems involving larger DOFs. To this aim, a large number of pseudo-observations are generated by fictitiously increasing the number of aircraft measurement. For this scenario, the domain covers California between 00:00Z and 23:59:59Z on 22 June 2008.  Four aircraft trajectories from the ARCTAS-CARB campaign are spread out over California as described in \citet{Guerrette17}.  A hypothetical plane is launched on each of those trajectories at approximately the beginning of each hour between 02Z and 15Z (9PM PST on 21 Jun 2008 to 10AM PST on 22 Jun 2008 PST) for a total of 14*4 = 56 flights.  This setup provides a total of 12,796 grid-scale observations of black carbon. As a result, the large BC problem has a higher DOFs of about 219 (see Fig. \ref{fig:bclarge_eigs}).

These observations were simulated in a base case that uses FINN biomass burning (BB) emissions and NEI 2005 anthropogenic (ANT) emissions. In our observation system simulation experiment (OSSE), the prior was defined by perturbing the inventories using a factor of 3.8 uncertainty for BB and factor of 2 for ANT, both with a 36 km horizontal error correlation length scale.  The BC observations were perturbed according to their uncertainties with no correlation. Two independent additive error sources are assumed for the observations: a uniform 12\% relative uncertainty and an absolute uncertainty of 0.02 $\mu$g/m$^3$. Other than the observations and OSSE setup, the experiment follows the settings for FINN$_\text{STD}$ described in \citet{Guerrette17}. As explained in Section \ref{smallBC}, the lognormal transformation applied to the emission control vector makes the cost function to minimize non-quadratic, which justifies the use of the incremental 4D-Var approach.

\begin{figure*}[]
\centering
 \includegraphics[width=43pc]{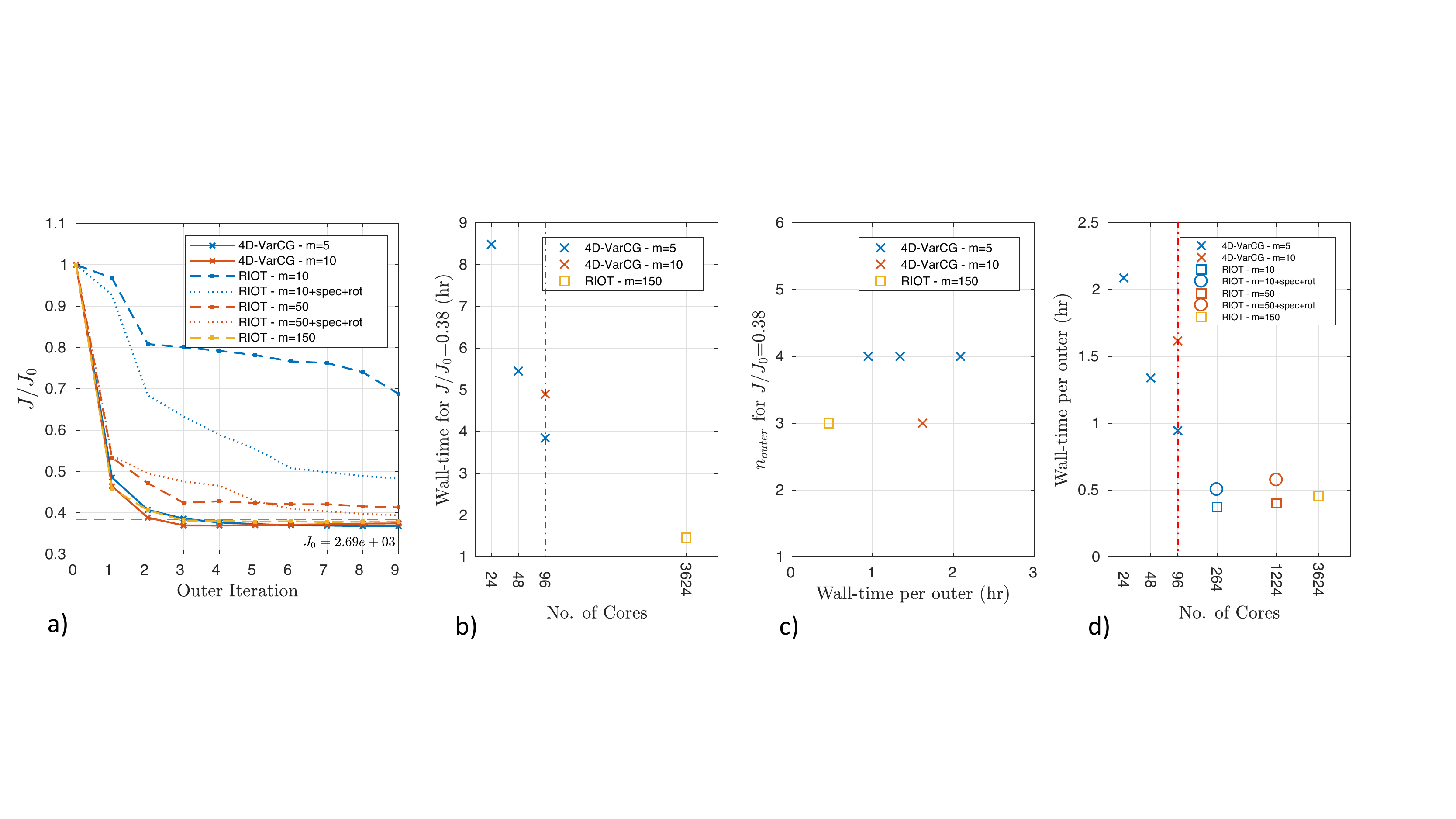}
\caption{Convergence performance for different configurations of RIOT and VarCG for the large BC problem. The VarCG experiments differ by the number of inner-loop iterations (m) and/or the number of cores used per model integration. The RIOT experiments differ by the rank of the Hessian approximation (m) and by the use of the spectral (spec) preconditioning and/or the rotation method (rot). a) Relative cost function value (fraction of initial value) (y-axis) as a function of the outer iteration number (x-axis); b) wall-time to reach a cost function reduction of 62\% (y-axis) and corresponding number of cores per model integration (x-axis); c) number of outer iterations to reach a cost function reduction of 62\% (y-axis) and corresponding wall-time (x-axis); d) wall-time per outer iteration (y-axis) and corresponding number of cores per model integration (x-axis). The vertical red dotted line shows the maximum number of cores above which the MPI parallelization of the WRF-DA system does not reduce the computational wall-time of the integration.}
\label{fig:bclarge}
\end{figure*}

\paragraph{Results and Discussion}
Figure \ref{fig:bclarge} shows the convergence performance for the non-quadratic cost function for different VarCG and RIOT configurations, as a function of the elapsed walltime of the minimization, the number of outer loop iterations or the number of cores. Similarly to previous results for the small BC experiment, the threshold for convergence was defined based on the maximum cost function reduction across all methods, such that all cost functions that are within $3$\% of the minimum cost function after 15 iterations are considered as converged. 

Both VarCG with $m=10$ inner iterations and RIOT with $m=150$ samples show the fastest convergence in 4 outer iterations, followed by VarCG with $m=5$ inner iterations in 3 outer iterations. The improvement from the preconditioning technique in RIOT is evident especially when $m=10$ samples are used in the RSVD. In the case with $m=50$ samples, the preconditioning allows further decrease of the cost function beyond outer iteration 5, in contrast with the non-preconditioned case for which the cost function plateaus. Preconditioning RIOT when $m=150$ samples are used in the RSVD degrades the performance of the minimization (not shown here), which is consistent with the results obtained for the small BC problem when the rank of the Hessian approximation was increased (see Sec. \ref{smallBC}). Again, this shows that the threshold to determine if the preconditioning should be used in RIOT or not for a given number of samples will depend on the specific inverse problem at hand. Figure \ref{fig:bclarge_eigs} shows the eigenspectrum of the prior-preconditioned Hessian in the first outer iteration of the minimization. Similar to what was found in the case of the small BC problem (see Sec. \ref{smallBC}), the Hessian eigenspectrum can help illustrate the impact of the preconditioner in RIOT and can be used as a diagnostic tool to adapt the preconditioning technique to specific inverse problems. In this example, the preconditioner is beneficial for ranks of the Hessian approximation corresponding to eigenvalue indices situated in the fast decreasing regime of the spectrum (e.g., $m=10$), while for higher ranks ($m=150$) the non-preconditioned minimization performs better.

Figure \ref{fig:bclarge} also shows the impact of exploiting the model-level parallelization in VarCG.  As the number of cores per model integration is increased from 24 to 96 (i.e., by a factor 4) in VarCG with $m=5$ inner iterations, the wall-time to reach convergence is reduced by more than a factor of 2 (see pannel b)). As discussed in Section \ref{smallBC}, the MPI parallelization limit for WRF-DA being defined by 96 cores, these results show 
that RIOT can be used to further reduce the wall-time cost by a factor 2.5 at the expense of mobilizing about 37 times more cores simultaneously during the RSVD parallel integrations. 

\begin{figure*}[]
\centering 
 \includegraphics[width=40pc]{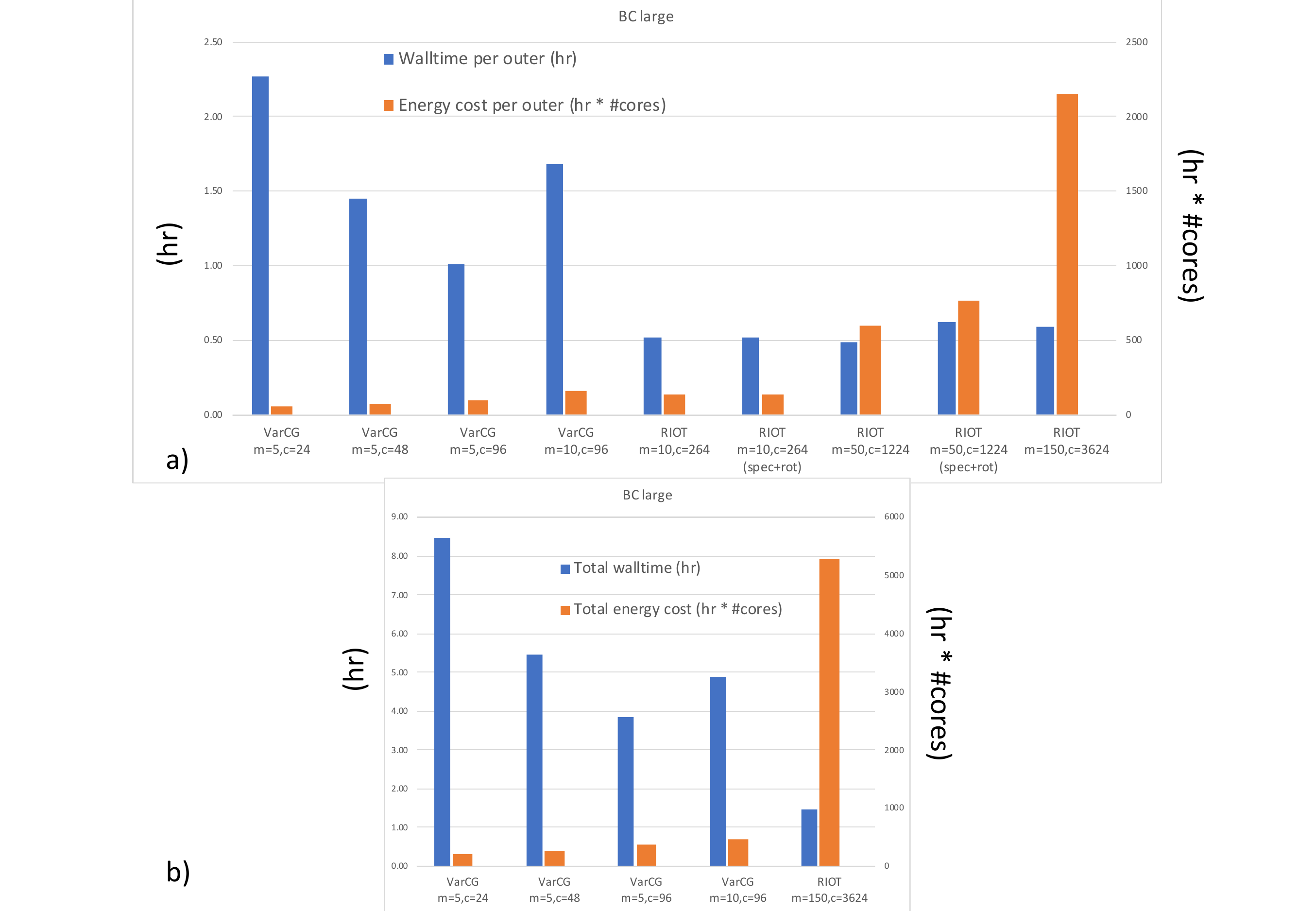}
\caption{Wall-time and energy cost performance for different configurations of RIOT and VarCG for the large problem. The number of inner iterations for VarCG or the rank of the Hessian approximation for RIOT are denoted by $m$, while the number of cores used per model integration is denoted by $c$. a) Wall-time (hours) and energy cost (hours$\times$cores) per outer iteration; b) wall-time (hours) and energy cost (hours$\times$cores) to reach a cost function reduction of 62\%. }
\label{fig:bclarge_chart}
\end{figure*}

Similar to the results shown in Sec. \ref{smallBC} for the small BC problem, Fig. \ref{fig:bclarge_chart} summarizes the performances of the different minimization configurations in term of both the wall-time and a metric related to energy consumption (i.e., the product of the number of cores and the elapsed walltime). Although RIOT with $m=10$ samples allows to drastically reduce the wall-time per outer iteration compared to VarCG without increasing the energy cost, a significant increase in the number of outer iterations is required to reach convergence when less than $m=150$ samples are used. As a consequence, at least $m=150$ samples are necessary to take advantage of the RIOT algorithm, in which case the wall-time can be reduced by a factor of about 2.5 compared to the best VarCG configuration ($m=5$ with $c=96$ cores), at the expense of about an order of magnitude increase in energy cost.

Note that in our experiment the cost of the TL and non-linear forward model integrations are approximately the same. It is worth mentioning that in an operational context (e.g., ECMWF's system) the TL and AD models in incremental 4D-Var usually correspond to a simplified version of the non-linear forward model and are run on a much coarser grid resolution to make the minimization computationally tractable. Therefore, similar numerical experiments conducted with such systems could lead to different conclusions. 



\section{Conclusion}
In this study a new technique for incremental 4D-Var systems, the Randomized Incremental Optimal Technique (RIOT), is introduced and evaluated using both a Lorenz-96 (L-96) system and a realistic black carbon source inversion problem based on WRF-DA. This method allows one to replace the standard iterative Conjugate-Gradient (VarCG) method in each quadratic inner loop minimization by a fully parallel randomized singular vector decomposition (RSVD) algorithm.
RIOT exploits the optimality results established by  \citet{bous18} and \citet{spantini2015} to compute the dominant modes of the prior-preconditioned Hessian using an RSVD method.  Those eigenmodes are used to build approximations of both the posterior increments and posterior error covariance matrices, which are adaptively constructed based on the characteristics of the Hessian spectrum.

For the L-96 problem, quasi-linear to very non-linear regimes were considered. Provided enough samples (i.e., parallel tangent-linear and adjoint integrations) are used in the RSVD, RIOT shows non-quadratic optimization performances comparable to the VarCG method. Our results demonstrate that the computational efficiency of RIOT relative to VarCG increases with the non-linearity of the L-96 system. The need for an outer loop Hessian preconditioning technique in RIOT has also been highlighted as well as the importance of adapting the traditional spectral preconditioner to the RSVD approach. In our experiments, the Hessian spectral preconditioner is combined with a so-called rotation technique in order to improve the RSVD sampling across outer iterations. For the L-96 problem, this preconditioning technique significantly improved the performance of the optimization and proved to be much more stable than either the spectral preconditioner alone or a non-preconditioned minimization, especially when the non-linearity of the system is increased.

For the black carbon source inversion problems, our analysis shows that for a given rank of the approximation (i.e., number of inner iterations in VarCG and number of RSVD samples in RIOT) the VarCG method generates posterior increments with systematically smaller errors than those obtained from RIOT, which stems from the minimum error residual property enforced by construction in VarCG.  However, RIOT is shown to produce equivalent solutions using a Hessian approximation of rank about three times as large as the number of VarCG inner iterations. Assuming a large number of RSVD samples can be run in parallel (i.e., a large amount of cores can be used simultaneously), the overall wall-time cost of the minimization in the black carbon inversion problems can be reduced by a factor 2 to 3 when using RIOT instead of an iterative VarCG minimization, at the expense of increasing the energy cost by a factor 4 to 10. Increasing the number of observations in the black carbon pseudo-inversion experiment so that the DOFs is at least an order of magnitude higher doubled the energy cost of the minimization. Our results also showed that the use of the spectral preconditioner with rotation can be detrimental to the minimization compared to a non-preconditioned approach when the rank of the Hessian approximation becomes large and the Hessian has a fast-decreasing spectra. Therefore, for a given inverse problem, a more tailored preconditioning approach based on the characteristics of its Hessian eigenspectra will be necessary to ensure efficient RIOT implementations. In the case of the posterior covariance errors estimates, RIOT is shown to be even more competitive with respect to the VarCG method, with about a factor of 40 improvement in wall-time performance for similar quality of the approximation. 

In order to assess the potential of RIOT for operational NWP DA applications, further numerical experiments are in progress based on ECMWF's incremental 4D-Var system. One limitation of RIOT is that the performance of the minimization will strongly depend on the characteristic of the Hessian spectra. Although inverse problems with fast decreasing Hessian spectra are associated with good performances of the RSVD approach, DA problems in NWP are known to generate flatter Hessian spectra, in particular due to the availability of large amount of data with similar weights in the assimilation. In such cases a compromise between Krylov subspace iterative algorithms and randomized parallel methods, such as the use of randomized power iterations \citep{Halko11} or the Block-Lanczos algorithm \citep{Golub89} to approximate the Hessian, is a promising approach to mitigate the dependence of the performance of the RIOT algorithm on the Hessian spectra. Preliminary results obtained from a Block-Lanczos approach and presented in this study in the context of the L-96 problem demonstrate the strong potential of such hybrid iterative-randomized methods and should be further investigated in the context of more realistic large-scale inverse problems. Another useful application of the RSVD approach in the incremental 4D-Var context would be to use randomized Hessian decompositions to precondition the outer loop. Indeed, the possibility to generate Hessian approximations of much higher rank using randomized decomposition techniques compared to iterative Krylov methods (inherently limited by their sequential nature) has a strong potential. Moreover, posterior square-root sampling methods used to allow for flow-dependency in cycling 4D-Var DA and that currently rely on Hessian approximations based on iterative minimization algorithms (e.g., BFGS or VarCG) are limited by the low-rank nature of the estimates, as suggested by \citet{Auligne16}. The performance of the RSVD to approximate the posterior error covariance matrix presented in this study clearly demonstrates the potential of this approach to address those limitations. 

\paragraph{Acknowledgments} 
 The authors acknowledge support from the NOAA grant NA16OAR4310113. The authors acknowledge the NOAA Research and Development High Performance Computing Program for providing computing and storage resources that have contributed to the research results reported within this paper (url:https://rdhpcs.noaa.gov). This research was performed while J. Guerrette held an NRC Research Associateship award at NOAA/ESRL. Nicolas Bousserez performed part of the research while working on the CHE project at ECMWF. The CHE project has received funding from the European Union's Horizon 2020 research and innovation programme under grant agreement No 776186.

\appendix
\section*{Appendix}
 \begin{algorithm}
  \caption{One Cycle of the Randomized Incremental Optimal Technique for Physical Space Statistical Analysis (RIOT-PSAS)}\label{alg:riot_dual}  
     \begin{algorithmic}[1]                    
        \STATE Start with $\mathbf{x}_{0}=\mathbf{x}_b$, $\mathbf{R}^{-1}=\mathbf{L}_{\mathbf{R}^{-1}}\mathbf{L}_{\mathbf{R}^{-1}}^\top$.
        \STATE Choose $m$ (number of samples for randomized SVD).
        \STATE Choose $k_f$ (number of outer loops).
        \FOR{$k=1, 2,\hdots, k_f$}         
         \STATE Integrate and store trajectory $H\left(\mathbf{x}_{k-1}\right)$  .               
          \STATE Compute and store $\mathbf{d}_{k-1}=\mathbf{y}-H(\mathbf{x}_{k-1})$.
           \STATE Let $\hat{\mathbf{A}} \equiv \mathbf{L}_{\mathbf{R}^{-1}}\mathbf{H}_{k-1}\mathbf{B}\mathbf{H}_{k-1}^\top \mathbf{L}_{\mathbf{R}^{-1}}^\top \in \mathbb{R}^{p \times p}$, where $\mathbf{H}_{k-1}$ is the tangent-linear and $\mathbf{H}_{k-1}^\top$ the adjoint model at $\mathbf{x}^{k-1}$.
        	   \STATE Draw $\mathbf{\Omega}_k\in\mathbb{R}^{p \times m}\sim\mathcal{N}(0,1)$.
       
	  \FORALL {$i\in\{1,2,\hdots,m\}$}  
	  \STATE $\mathbf{\omega}_{i}=\mathbf{\Omega}_k(:,i)$.
	   \ENDFOR
           \FORALLP {$i\in\{1,2,\hdots,m\}$}
                
                 \STATE $\mathbf{y}_{i} = \hat{\mathbf{A}}\mathbf{\omega}_{i}$   
                  \STATE Compute $\mathbf{b}'\equiv \mathbf{H}_{k-1}( \mathbf{x}_b-\mathbf{x}_{k-1})$.

           \ENDFOR           
           \STATE Form $\mathbf{Y}=\left[\mathbf{w}_1,\mathbf{w}_2,\hdots,\mathbf{w}_{m}\right]$.
           \STATE Compute $\mathbf{Q}\in\mathbb{R}^{p \times m}$ orthonormal basis of $\mathbf{Y}$ from, e.g., $\text{QR}$ algorithm.
           \STATE Solve for $\mathbf{K}$ in $\mathbf{K}\mathbf{Q}^\top\mathbf{\Omega} = \mathbf{Q}^\top\mathbf{Y}$.
           \STATE Form eigendecomposition, $\mathbf{K}=\mathbf{Z}\mathbf{\Lambda}\mathbf{Z}^\top$ ($\mathbf{\Lambda} = diag(\lambda_i)$).
                      \STATE Compute  $\mathbf{U}=\mathbf{Q}\mathbf{Z}$.

  \IF{ $\lambda_m>1$ }
         \STATE $\delta\mathbf{w}=  \mathbf{L}_{\mathbf{R}^{-1}} \left (\sum_{i=1}^m \frac{1}{1+\lambda_i}\mathbf{u}_i\mathbf{u}_i^\top \right )\mathbf{L}_{\mathbf{R}^{-1}}^\top\left ( \mathbf{d}_{k-1}-\mathbf{b}' \right ) $
  \ELSE
                \STATE $\delta\mathbf{w}=  \mathbf{L}_{\mathbf{R}^{-1}} \left (\mathbf{I}-\sum_{i=1}^m \frac{\lambda_i}{1+\lambda_i}\mathbf{u}_i\mathbf{u}_i^\top \right )\mathbf{L}_{\mathbf{R}^{-1}}^\top\left ( \mathbf{d}_{k-1}-\mathbf{b}' \right ) $
  \ENDIF
          \STATE $\mathbf{x}_{k} = \mathbf{x}_{k-1}+\mathbf{B}\mathbf{H}_{k-1}^\top\delta\mathbf{w}$
        \ENDFOR
    \end{algorithmic}
    \end{algorithm}


\newpage
\bibliography{ms}{}
\bibliographystyle{wileyqj}

\end{document}